\documentclass{article}

\usepackage{indentfirst,csquotes}

\usepackage{graphicx} 
\usepackage{amsmath,amssymb,amsfonts,amsthm}
\usepackage{hyperref}
\usepackage{xcolor}
\usepackage{comment}
\usepackage[sortcites=true]{biblatex}
\theoremstyle{definition}
\newtheorem{definition}{Definition}[section]
\newtheorem{setting}[definition]{Setting}
\newtheorem{remark}[definition]{Remark}
\theoremstyle{plain}
\newtheorem{theorem}[definition]{Theorem}
\newtheorem{lemma}[definition]{Lemma}
\newtheorem{corollary}[definition]{Corollary}
\newtheorem{claim}[definition]{Claim}

\newcommand{\R}{\mathbb{R}}
\newcommand{\Z}{\mathbb{Z}}
\newcommand{\N}{\mathbb{N}}
\newcommand{\C}{\mathbb{C}}
\newcommand{\Q}{\mathbb{Q}}

\providecommand{\keywords}[1]{\small\textbf{\textit{Keywords:}} #1}

\title{Full Rank and Stable Sampling of~Tanh Networks}
\date{\today}

\author{%
  Philipp Grohs\textsuperscript{1,2} and Matěj Trödler\textsuperscript{1}\\[0.5em]
  {\small \textsuperscript{1}Faculty of~Mathematics, University of~Vienna}\\
  {\small \textsuperscript{2}RICAM, Austrian Academy of~Sciences}
}

\begin{document}

\maketitle

\begin{abstract}
    We study the recovery of feedforward $\tanh$ networks of fixed architecture with $P$ parameters from finitely many point evaluations. First, we characterize the set of parameters where the differential of the realization map has full rank and show that it coincides with the regular parameters in the sense of Vla{\v{c}}i{\'c} and B{\"o}lcskei (2021). 
    Second, we make use of this characterization to prove that almost every tuple of $2P+1$ sampling points yields uniform-norm Lipschitz stability on compact subsets of this parameter set. Third, we prove that for exact observations, projected gradient descent and projected Gauss--Newton converge locally from sufficiently close feasible initial guesses. The same evaluations also determine the realized function for arbitrary parameters and identify parameters in the characterized set up to permutations and sign changes of hidden neurons. Finally, we develop a~general framework for real-analytic networks that highlights the role of full rank in stability and local reconstruction.
\end{abstract}

\keywords{tanh neural networks, full rank, stable sampling, finite-sample identifiability, projected gradient descent (PGD), projected Gauss--Newton (PGN)}

\section{Introduction}

Deep neural networks provide efficient approximations for a~range of~function classes \cite{elbrachter2021deep,berner2021modern}. 
In practice, however, the information available about a function to be approximated often consists of only finitely many samples. A central problem is therefore to determine an accurate neural network approximation from these samples. In this paper, we establish a link between stable approximation from samples and the degeneracy of the derivative of the neural network realization map with respect to its parameters. Motivated by this connection, we first characterize the parameters at which this derivative has full rank for $\tanh$ networks of fixed architecture with $P$ parameters and show that they are precisely the regular parameters in the sense of \cite{vlavcic2021affine}. Using this characterization together with the analytic regularity of $\tanh$ neural networks we then prove that, on any compact subset $\Xi$ of this regular parameter set, for almost every tuple $(x_i)_{i=1}^{2P+1}$ of $2P+1$ sampling points, a norm equivalence of the form $$\|u_\theta - u_\nu\|_{L^\infty}\sim \max_{1\le i\le 2P+1}\|u_\theta(x_i)-u_\nu(x_i)\|_2\quad\mbox{for all}\quad \theta,\mu \in \Xi$$ holds true with $u_\theta$ denoting the $\tanh$ neural network with parameters $\theta\in \Xi$. Using this result we can show that 
a quasi-best uniform approximation of any continuous function in the class $(u_\theta)_{\theta\in \Xi}$ can be computed from $2P+1$ generic samples. To our knowledge this is the first (quasi-) optimal sampling result for classes of deep neural networks. 
Finally, our analysis also allows to study the algorithmic reconstruction of $\tanh$ neural networks from samples using empirical risk minimization.

\subsection{Main results}

We now present our main results for feedforward $\tanh$ networks of~arbitrary depth and width, with an~affine output layer. We show that the differential of~the realization map has full rank at almost every parameter, in contrast to the rank deficiency caused by rescaling symmetries in ReLU networks. This result underpins our stability and local reconstruction results: a~number of~point evaluations linear in the parameter count suffices to control the uniform reconstruction error on compact parameter sets where the differential has full rank and to ensure local convergence of~projected gradient descent (PGD) and projected Gauss--Newton (PGN). The same evaluations also determine the realized function and identify parameters satisfying the full rank condition up to permutations and sign changes of~hidden neurons. We first specify the architecture, parameter symmetries, and regular parameters. Throughout, $[m]:=\{1,\dots,m\}$, and for compact $K\subset\R^d$ and $f\in C(K, \R^n)$, we write
\[
    \|f\|_{L^\infty(K)}:=\sup_{x\in K}\|f(x)\|_2.
\]
We denote the Hausdorff dimension of~a~set $\Theta$ by~$\dim_H(\Theta)$.

\begin{setting}\label{set:deep_tanh_nets}
    \noindent\textbf{Architecture.} Let $d,L\in\N$, set $n_0:=d$, and let $n_\ell\in\N$ for $\ell\in[L+1]$. For $\ell\in[L+1]$, let
    \[
        W^\ell=(w_{ij}^\ell)_{i,j}
        \in\R^{n_\ell\times n_{\ell-1}},
        \qquad
        b^\ell\in\R^{n_\ell}.
    \]
    Set $h^0(x):=x$ and define recursively
    \begin{equation}\label{eq:deep_tanh_layers}
        \begin{aligned}
            q^\ell(x)
            &:=b^\ell+W^\ell h^{\ell-1}(x),
            &&\ell\in[L+1],\\
            h^\ell(x)
            &:=\tanh(q^\ell(x)),
            &&\ell\in[L],
        \end{aligned}
    \end{equation}
    where $\tanh$ is applied componentwise. Set $n:=n_{L+1}$ and $P:=\sum_{\ell=1}^{L+1} n_\ell(n_{\ell-1}+1)$. For $\theta=(W^\ell,b^\ell)_{\ell=1}^{L+1} \cong \R^P$, the realized function $u_\theta\colon\R^d\to\R^n$ is defined by
    \[
        u_\theta(x)=\mathcal R^{\tanh}(\theta)(x):=q^{L+1}(x).
    \]

    \noindent\textbf{Parameter symmetries.} For every $\ell\in[L]$, let $P_\ell\in\R^{n_\ell\times n_\ell}$ be a~permutation matrix. Let
    \[
        D_\ell
        :=
        \operatorname{diag}
        (\varepsilon_1^\ell,\dots,\varepsilon_{n_\ell}^\ell),
        \qquad
        \varepsilon_i^\ell\in\{\pm1\},
    \]
    and set
    \[
        Q_0:=I_d,
        \qquad
        Q_\ell:=D_\ell P_\ell,\quad \ell\in[L],
        \qquad
        Q_{L+1}:=I_{n_{L+1}}.
    \]
    Finally, let $\Gamma^{\tanh}$ be the set of~all parameter transformations
    obtained in this way. For $\theta'=\gamma(\theta)$, $\gamma\in\Gamma^{\tanh}$, the transformed parameters are given by
    \[
        (W^\ell)'
        =
        Q_\ell W^\ell Q_{\ell-1}^{T},
        \qquad
        (b^\ell)'
        =
        Q_\ell b^\ell,
        \qquad
        \ell\in[L+1].
    \]
\end{setting}

The transformations in $\Gamma^{\tanh}$ form a~finite group of~diffeomorphisms and preserve the realized function:
\[
    u_{\gamma(\theta)}=u_\theta
    \qquad\text{for all }\theta\in\R^P
    \text{ and }\gamma\in\Gamma^{\tanh}.
\]
Indeed, permuting hidden neurons preserves the output, and the oddness of~$\tanh$ allows us to change the sign of~a~neuron's incoming weights and bias together with its outgoing weights. We write $\mu\sim_{\Gamma^{\tanh}}\nu$ if $\mu=\gamma(\nu)$ for some $\gamma\in\Gamma^{\tanh}$. We therefore study parameter recovery on the regular set defined below, modulo these symmetries.

\begin{definition}[Regular parameters]\label{def:regular_parameters}
    Let $\Theta_{\mathrm{reg}}^{\tanh}\subset\R^P$ denote the set of~parameters satisfying the following conditions:
    \begin{enumerate}
        \item for every $\ell\in[L]$ and $i\in[n_\ell]$,
        \begin{align}
            W^\ell_{i,:}&\neq0,
            \label{eq:every_neuron_has_parent}\\
            W^{\ell+1}_{:,i}&\neq0;
            \label{eq:every_neuron_has_outgoing_edge}
        \end{align}
    
        \item for every $\ell\in[L]$ and every distinct
        $i,j\in[n_\ell]$,
        \begin{equation}\label{eq:regularity_of_params_2}
            \bigl(W^\ell_{i,:},b_i^\ell\bigr)
            \neq
            \pm
            \bigl(W^\ell_{j,:},b_j^\ell\bigr).
        \end{equation}
    \end{enumerate}
\end{definition}

These conditions correspond to regularity in the sense of~\cite[Definition~8]{vlavcic2021affine}. Condition \eqref{eq:every_neuron_has_parent} ensures that the designated input nodes are the only nodes without incoming edges. Condition \eqref{eq:every_neuron_has_outgoing_edge} ensures non-degeneracy by~giving each hidden neuron a~path to the output. Condition \eqref{eq:regularity_of_params_2} ensures $\tanh$-irreducibility by~excluding pairs of~neurons in the same layer whose incoming weights and biases agree up to a~sign. Note that the regular parameter set is the complement of~a~finite union of~proper linear subspaces of~$\R^P$. It is therefore open, dense, and of~full Lebesgue measure, with $\dim_H(\Theta_{\mathrm{reg}}^{\tanh})=P$.

We now present our three main results on full rank of the $\tanh$ realization map, stable sampling, and convergence of PGD and PGN. We also state Corollary~\ref{col:finite_sample_identifiability_tanh}, which follows directly from known results but, to the best of our knowledge, has not previously been stated explicitly.

\subsubsection{Full rank of~the realization map}

Our main technical result is that the differential of~the realization map $\theta\mapsto u_\theta$ has full rank at every regular parameter of~a~deep $\tanh$ network. Section~\ref{sec:gen_results} illustrates the role of~full rank by~showing how it is used to establish stable sampling and local convergence of~reconstruction algorithms. The following theorem verifies this assumption for regular deep $\tanh$ networks.

\begin{theorem}[Full rank]\label{thm:full_rank_realization_d}
    It holds that
    \[
        \Theta^{\tanh}_{\mathrm{reg}} = \left\{
        \theta\in\R^P:
        \left\{
        \frac{\partial u_\theta}{\partial\theta_j}
        \right\}_{j\in[P]}
        \text{ is linearly independent in }C(\R^d, \R^n)
        \right\}.
    \]
    Equivalently, the range of $D\mathcal R^{\tanh}(\theta)$ has dimension $P$ if and only if $\theta\in\Theta^{\tanh}_{\mathrm{reg}}$.
\end{theorem}

As in the ReLU example above, each failure of regularity yields a nonzero direction in the kernel of $D\mathcal R^{\tanh}(\theta)$ through a parameter variation that preserves the realized function:
\begin{itemize}
    \item If \eqref{eq:every_neuron_has_parent} fails, a~neuron has zero incoming weights and therefore a~constant output. Changes in its outgoing weights can then be absorbed into the next layer's bias.
    \item If \eqref{eq:every_neuron_has_outgoing_edge} fails, a~neuron has zero outgoing weights, so its bias can vary freely.
    \item If \eqref{eq:regularity_of_params_2} fails, two neurons in the same layer have outputs that agree up to sign, so their outgoing weights can be varied in compensating directions.
\end{itemize}    
Hence, each of~these variations gives a~nonzero direction in the kernel of~$D\mathcal R^{\tanh}(\theta)$, so the differential cannot have full rank at the corresponding parameters. The main difficulty is proving the converse, namely that the parameter derivatives are linearly independent at every regular parameter. Establishing this implication is one of~the main contributions of~this work.

Moreover, combining the full-rank theorem above with Theorem~\ref{thm:generic_local_sampling}, whose proof uses the finite witness theorem of~\textcite{amir2023neural}, we obtain that for Lebesgue almost every sampling tuple $X=(x_1,\ldots,x_{2P+1})$, the sampled Jacobian
\[
    J(\theta)=D_\theta\bigl(u_\theta(x_1),\ldots,u_\theta(x_{2P+1})\bigr)\in\R^{n(2P+1)\times P}
\]
has full column rank simultaneously at all regular parameters.

The range of~this differential is the tangent space generated by~parameter variations. Its full rank is a~key ingredient in our stable-sampling and local reconstruction results. Related questions arise in neural Galerkin methods for time-dependent PDEs, where the dynamics are projected onto the span of~the parameter derivatives \cite{bruna2024neural}. In this setting, \cite{zhang2026sequential} investigates "tangent space collapse", including linear dependence and numerical near-dependence of~these derivatives, and its consequences for accuracy.

\subsubsection{Stable sampling}

The result above is the main ingredient for our following result, which is stability of~learning these networks. In particular, on compact sets of~regular parameters, we show that the $L^\infty(K)$ difference of~two networks can be controlled by~the differences at the sampled points. The bound holds for almost every choice of~$2P+1$ sample points, so the number of~samples grows only linearly in $P$, which is particularly interesting.

\begin{theorem}[Stable sampling]
\label{thm:stable_sampling_tanh}
Let $K\subset\R^d$ be compact, and let $W\subset\operatorname{int}(K)$ be nonempty and open. For Lebesgue almost every $X=(x_1,\dots,x_{2P+1})\in W^{2P+1}$, the following holds: for every nonempty compact set $\Xi\subset\Theta_{\mathrm{reg}}^{\tanh}$, there is a~constant $\delta_{\Xi,X}>0$ such that
\[
    \max_{1\le i\le 2P+1} \|u_\mu(x_i)-u_\nu(x_i)\|_2
    \ge \delta_{\Xi,X}
    \|u_\mu-u_\nu\|_{L^\infty(K)}
    \qquad \text{for } \mu,\nu\in\Xi.
\]
\end{theorem}

This result fits naturally into the theory of~sampling discretization, which studies how function norms can be controlled by~finitely many point evaluations. Uniform-norm discretization for finite-dimensional linear spaces is investigated, for instance, in \cite{kashin2022sampling,krieg2026sampling}. Our theorem establishes an~analogous inequality for the nonlinear difference class $\mathcal D_\Xi:=\{u_\mu-u_\nu:\mu,\nu\in\Xi\}$. The proof combines differential nondegeneracy, which controls nearby parameters, with finite-sample identifiability and compactness. However, the number of~samples needed to achieve small error may be still large. For instance, for the ReLU approximation spaces studied in \cite{grohs2024proof}, guaranteeing a~worst-case uniform error of~at most $\varepsilon$ requires at least on the order of~$\varepsilon^{-d}$ samples. Nevertheless, this result is particularly interesting that only linear number of~samples can provide the concrete approximation guarantees.

The estimate also provides a~stability guarantee for reconstruction from noisy observations. Suppose $\theta_\star,\widehat\theta\in\Xi$ and, at every sampling location, the observation and fitting errors satisfy $\max_{i\in[2P+1]} \| y_i-u_{\theta_\star}(x_i) \|_2 \le\varepsilon$ and $\max_{i\in[2P+1]}  \|u_{\widehat\theta}(x_i)-y_i \|_2 \le\eta$, respectively. The triangle inequality and the stable-sampling estimate give
\[
    \|u_{\widehat\theta}-u_{\theta_\star}\|_{L^\infty(K)}
    \le \frac{\varepsilon+\eta}{\delta_{\Xi,X}}.
\]
Thus, $\delta_{\Xi,X}^{-1}$ controls the amplification of~observation and fitting errors when estimating the error over the full domain.

Finally we mention that Theorem \ref{thm:stable_sampling_tanh} allows to determine a quasi-best uniform approximation of any $f\in C(K)$ from $2P+1$ generic samples: For $f\in C(K)$ denote
$$
    \theta_\ast\in \mathop{\mbox{argmin}}_{\theta\in \Xi}\|f-u_\theta\|_{L^\infty(K)},\quad 
    \hat\theta\in \mathop{\mbox{argmin}}_{\theta\in \Xi}\max_{1\le i\le 2P+1}\|f(x_i)-u_\theta(x_i)\|_{2}.
$$
Then, Theorem \ref{thm:stable_sampling_tanh} together with the triangle inequality readily yields that
$$
    \|f- u_{\hat\theta}\|_{L^\infty(K)}\le \left(1+\frac{2}{\delta_{\Xi,X}}\right)\|f - u_{\theta_\ast}\|_{L^\infty(K)}.
$$
In other words, \emph{a quasi-best uniform approximation of $f$ can be constructed from $2P+1$ generic samples.}
\subsubsection{Locally convergent reconstruction}

The preceding stable-sampling estimate bounds the uniform reconstruction error by~the sum of~the observation and fitting errors multiplied by~$\delta_{\Xi,X}^{-1}$, which can be large. This raises the question of~whether standard algorithms can find such a~fit, given the known difficulties of~recovering neural networks from samples \cite{grohs2024proof,grohs2025theory,berner2023learning,grohs2026limitations}. The least-squares loss is generally nonconvex. At an~exact fit, its Hessian is positive definite whenever the sampled Jacobian has full column rank. Using this property, we prove that projected gradient descent and projected Gauss--Newton converge locally from sufficiently close feasible initial guesses when the observations are exact.

For $X=(x_1,\dots,x_{2P+1})$, define the sampling map and least-squares loss by
\[
    F_X(\theta):=(u_\theta(x_i))_{i=1}^{2P+1}\in\R^{n(2P+1)},
    \qquad
    L_{\theta_\star}(\theta)
    :=\frac12\|F_X(\theta)-F_X(\theta_\star)\|_2^2.
\]
Moreover, define $J(\theta):=DF_X(\theta)\in\R^{n(2P+1)\times P}$. For a~nonempty closed convex set $C\subset\R^P$, let $P_C$ denote Euclidean projection onto $C$.

\begin{theorem}[Local convergence of~PGD and PGN]
\label{thm:local_reconstruction_tanh_simple}
Let $W\subset\R^d$ be nonempty and open. For Lebesgue almost every $X\in W^{2P+1}$, the following holds. For every nonempty compact set $\Xi\subset\Theta^{\tanh}_{\mathrm{reg}}$, there exist $r,\alpha_0,\beta>0$ such that, for every $\theta_\star\in\Xi$, every nonempty closed convex set $C\subset\mathbb R^P$, every $\vartheta\in C\cap\Gamma(\theta_\star)$, every fixed step size $\alpha\in(0,\alpha_0]$, and every $\theta_1\in C$ with $\|\theta_1-\vartheta\|_2<r$, the following holds:
\begin{itemize}
    \item \textbf{projected gradient descent:} the iteration
\[
    \theta_{k+1}
    =P_C\bigl(\theta_k-\alpha\nabla L_{\theta_\star}(\theta_k)\bigr)
\]
satisfies
\[
    \|\theta_k-\vartheta\|_2
    \le q_\alpha^{k-1}\|\theta_1-\vartheta\|_2,
    \qquad k\ge1,
    \qquad q_\alpha:=1-\alpha\beta\in(0,1).
\]
\item \textbf{projected Gauss-Newton:} the iteration
\[
    \theta_{k+1}
    =P_C\left(\theta_k-
    \bigl(J(\theta_k)^\top J(\theta_k)\bigr)^{-1}
    \nabla L_{\theta_\star}(\theta_k)\right).
\]
is well defined and converges to $\vartheta$, either in finitely many steps or at a~$Q$-superlinear rate.
\end{itemize}
\end{theorem}

The full PGD result also controls the loss and the uniform function error, with local constants uniform over targets in compact regular parameter sets; see Corollary~\ref{cor:local_pgd_c1}.

\subsubsection{Finite-sample identifiability of~tanh networks}

Combining the functional identifiability result of~\cite{vlavcic2021affine} with the finite witness theorem of~\cite{amir2023neural} yields the following corollary; see Lemmas~\ref{lem:global_identifiability_tanh} and~\ref{lem:generic_global_separation}. A single generic sampling tuple determines every realized function in the network class and, for regular parameters, determines the parameters up to $\tanh$-isomorphisms.

\begin{corollary}[Finite-sample identifiability]
\label{col:finite_sample_identifiability_tanh}
Let $W\subset\R^d$ be nonempty and open. For Lebesgue almost every $X=(x_1,\dots,x_{2P+1})\in W^{2P+1}$, the equivalence
\[
    u_\mu(x_i)=u_\nu(x_i)
    \text{ for every }i\in[2P+1]
    \quad\iff\quad
    u_\mu=u_\nu
    \text{ in }C(\R^d,\R^n)
\]
holds simultaneously for all $\mu,\nu\in\R^P$. If additionally $\mu,\nu\in\Theta_{\mathrm{reg}}^{\tanh}$, these conditions are equivalent to $\mu\sim_{\Gamma^{\tanh}}\nu$.
\end{corollary}

These conclusions concern exact point evaluations. In our earlier work \cite{grohs2026limitations}, we studied recovery of~$\tanh$ networks from quantized point evaluations. Under suitable assumptions on the architecture, parameter bounds, and quantization threshold, we proved that the worst-case expected $L^p([0,1]^d)$ reconstruction error of~any adaptive randomized algorithm using $m$ samples is bounded below by~a multiple of~$m^{-1/p}$, for sample budgets up to a~threshold that can grow exponentially with network depth. The present result complements these limitations by~establishing finite-sample identifiability of~regular networks from exact observations.

The same argument applies to networks with other real-analytic activations. By Lemma~\ref{lem:generic_global_separation}, almost every tuple of~$2P+1$ sampling points separates the realized functions. To obtain finite-sample parameter identifiability up to the corresponding symmetries, it therefore suffices to prove identifiability up to those symmetries on the whole domain.

\subsection{Related work}\label{sec:related_work}

We show that the differential of the realization map has full rank $P$ at every regular $\tanh$ parameter. For ReLU networks, the functional dimension is at most $P-\sum_{\ell=1}^{L}n_\ell<P$ \cite{grigsby2025functional}. Rank tests also play a~role in the study of local identifiability for deep ReLU networks \cite{bona2022local}. Related Hessian rank bounds for deep linear networks and certain shallow ReLU networks appear in \cite{singh2021analytic}. In our setting, generic sampling makes the least-squares Hessian positive definite at every regular exact fit, giving local strong convexity.

For deep $\tanh$ networks with $P$ parameters, almost every tuple of $2P+1$ sampling points yields uniform-norm Lipschitz stability on compact sets of regular parameters. However, this phenomena is not completely new. Using Christoffel sampling, \cite{cardenas2023cs4ml} obtain an~$L^2$ analogue with high probability and a~log-linear sample count. Moreover, \cite{alberti2022infinite} show that full-data Lipschitz stability persists under suitable finite-dimensional discretizations.

However, even when neural networks are identifiable, a~more practical question is whether they can be recovered efficiently from finitely many evaluations. Quantitative approximation bounds are available for $\tanh$ networks \cite{de2021approximation}, but strong approximation guarantees do not necessarily yield comparable recovery guarantees from samples. For the ReLU approximation class studied in \cite{grohs2024proof}, any algorithm that guarantees a~worst-case uniform error of~at most $\varepsilon$ requires $m\gtrsim\varepsilon^{-d}$ point evaluations. For $L^p$ recovery of ReLU networks, \cite{berner2023learning} establish a~worst-case error lower bound of order $m^{-1/d-1/p}$, even for adaptive randomized algorithms. In particular, uniform recovery gain requires $m\gtrsim\varepsilon^{-d}$. Further lower bounds for ReLU and ReQU approximation spaces appear in \cite{grohs2025theory,abdeljawad2025sampling}. For $\tanh$ networks observed through quantized point evaluations, \textcite{grohs2026limitations} establish an $m^{-1/p}$ error lower bound for sample budgets up to a~threshold that can grow exponentially with network size as well.

Several works also study the recovery of~shallow neural networks from data using gradient descent and its stochastic or projected variants, under different assumptions on the activation function, the samples, and the noise; see, for instance, \cite{zhong2017recovery,zhang2019learning,fornasier2025efficient,li2017convergence}. For deep networks, \cite{du2019gradient} show that gradient descent drives the training loss to zero at a~linear rate for sufficiently wide fully connected networks with random initialization. \cite{chatterjee2026convergence} proves convergence to an~interpolating solution for deep networks with smooth, strictly increasing activations. Stochastic Gauss--Newton has also been studied for training overparameterized deep networks \cite{cayci2026nonasymptotic}. 

For the projected methods, \cite{necoara2019linear} establish linear convergence of~projected gradient descent for convex problems under a~quadratic growth condition. \cite{izmailov2026local} discuss local superlinear convergence of~projected Gauss--Newton for constrained nonlinear equations when the Jacobian at the solution has full column rank. In our setting, a~single generic tuple of~$2P+1$ exact point evaluations gives a~sampled Jacobian with full column rank at every regular $\tanh$ parameter. We use this to prove that both projected methods converge from sufficiently close feasible initial guesses.

Finally, a~classical question in the theory of neural networks is whether the realized function determines the parameters up to symmetry. Let $\Gamma$ be a~group acting on $\Theta$ by parameter transformations that preserve the realized function. Given an observation set $\mathcal X\subset\R^d$, we ask under which assumptions the following implication holds for all $\mu,\nu\in\Theta$:
\[
    u_\mu(x)=u_\nu(x)
    \quad\text{for all }x\in\mathcal X
    \qquad \stackrel{?}{\Longrightarrow} \qquad
    \mu=\gamma(\nu)
    \quad\text{for some }\gamma\in\Gamma.
\]
A central issue is the choice of~the observation set $\mathcal X$: can it be finite, and how does its required cardinality depend on the network architecture?

Early work in the 1990s primarily addressed functional identifiability, assuming knowledge of~the entire input--output map, i.e., $\mathcal{X}=\R^d$. For shallow networks, that is, $L=1$, classical results for $\tanh$ include \cite{albertini1993uniqueness,sussmann1992uniqueness}. Related results for more general asymptotically constant activation functions were developed in \cite{kuurkova1994functionally}, while Gaussian and certain non-polynomial rational activations were studied in \cite{kainen1994uniqueness}. Identifiability and parameter symmetries of~ReLU networks have received considerable attention in recent years. For shallow networks, for instance, \cite{petzka2020notes} characterizes transformations that preserve the realized function.

For deep $\tanh$ networks, \cite{fefferman1994reconstructing} established identifiability under genericity assumptions on the parameters. These assumptions were substantially relaxed in \cite{vlavcic2021affine}, which establishes identifiability of~regular $\tanh$ networks modulo permutations and sign changes. Related identifiability results for a~family of~sigmoidal nonlinearities appear in \cite{vlavcic2022neural}. For deep ReLU networks, functional identifiability and query-based reconstruction have been investigated, for example, in \cite{bona2023parameter,rolnick2020reverse}.

Finite-sample identification has also been studied constructively under additional assumptions on the networks and sampling procedures; see, for example, \cite{fornasier2021robust,fornasier2025efficient}. For irreducible shallow ReLU networks, \cite{xia2025finite} shows that no universal finite sampling set separates all realized functions in the class, while also establishing target-dependent finite-sample identification under additional assumptions. For $\tanh$ and the logistic sigmoid, the same work constructs a~universal finite sampling set for irreducible shallow networks of~fixed width, although its sample bound is exponential in the hidden-layer width.

For deep ReLU networks, \cite{stock2023embedding} derives conditions for local identifiability from finitely many samples, while \cite{bona2022local} provides a~geometric characterization of~local identifiability and associated rank-based tests. For smooth deep networks with a~pyramidal architecture, \cite{fiedler2023stable} studies stable finite-sample recovery of~entangled weights. A related result in analytic inverse problems is due to \textcite{alberti2026analytic}, who show that, under a~separation assumption, $2r+1$ random scalar measurements almost surely identify each element of an $r$-dimensional analytic model. For further references on identifiability and stable sampling in this setting, see \cite[Section~1.4]{alberti2026analytic}.

\subsection{Outline}

Section~\ref{sec:gen_results} develops a~general framework for identifying realized functions from finitely many samples, establishing stable sampling, and proving local convergence of~projected gradient descent and projected Gauss--Newton. Section~\ref{sec:tanh_results} specializes these results to $\tanh$ networks. We first prove the private-poles lemma (Lemma~\ref{lem:private_poles}) and use it to establish full rank of~the realization map at every regular parameter (Theorem~\ref{thm:full_rank_realization_d}).

\section{General results}\label{sec:gen_results}

We now develop a~general framework for stable sampling and locally convergent reconstruction of~parameterized function families. The basic setting assumes local Lipschitz continuity of~the realization map and is not tied to a~particular activation function. The sampling results require additional regularity and subanalyticity assumptions, stated explicitly below. Real-analytic activations, such as the logistic sigmoid, softplus, Gaussian, GELU, and SiLU, provide natural examples for which the required analyticity in the input variable holds. Application of~the results to these networks additionally requires verification of~the remaining hypotheses, including full rank and, for global stable sampling, identifiability modulo the prescribed finite symmetry group. The $\tanh$ networks of~Setting~\ref{set:deep_tanh_nets} provide a~concrete instance of~this framework. 

We first introduce the general settings and definitions. Recall that for a~matrix $A\in\R^{m\times P}$, the induced Euclidean norm is given by
\[
    \|A\|_2:=\sup_{\|v\|_2=1}\|Av\|_2
    =\sigma_{\max}(A)
    =\sqrt{\lambda_{\max}(A^\top A)}.
\]
\begin{setting}\label{set:parmodel}
    Let $P,d\in \mathbb{N}$, $\Theta\subset \mathbb{R}^P$, and $K\subset \mathbb{R}^d$ compact. Consider a~realization map
    \begin{equation*}
        \mathcal{R}:\left\{\begin{array}{ccc}\Theta\subset \mathbb{R}^P& \to & C(K, \R^n)\\ \theta & \mapsto & u_\theta:=\mathcal{R}(\theta)\end{array}\right. .
    \end{equation*}
    For $\theta\in \Theta$ and $r\in (0,\infty)$ we denote $B(\theta,r):=\{\mu\in \Theta:\ |\theta - \mu|<r\}$.  
    We assume that $\mathcal{R}: \Theta \to C(K, \R^n)$ is locally Lipschitz: for every $\theta\in \Theta$, there exist $r, L>0$ such that 
    \begin{equation}\label{eq:lipschitz}
        \mbox{for every }\mu,\nu\in B(\theta,r):\ \|u_\mu - u_\nu\|_{L^\infty(K)}\le L |\nu-\mu|.
    \end{equation}
\end{setting}
\begin{definition}[Parameter symmetries, identifiability]
    Let $\Gamma \subset C^1(\Theta,\Theta)$ be a~finite group of~$C^1$ diffeomorphisms, leaving $\mathcal{R}$ invariant in the sense that
    \begin{equation}\label{eq:gamma_invariant}
        \mbox{for all }\theta\in \Theta,\ \gamma\in \Gamma:\ u_\theta = u_{\gamma(\theta)} 
    \end{equation}
    and write
    \[
        \theta \sim_{\Gamma}\mu \iff \exists \gamma\in \Gamma:\ \mu = \gamma(\theta). 
    \]
    We say that global identifiability holds for $(\mathcal{R},\Theta)$ modulo $\Gamma$ if for all $\nu,\mu\in \Theta$ it holds that
    \[
        \mathcal{R}(\mu) \equiv \mathcal{R}(\nu)
        \implies
        \mu \sim_{\Gamma} \nu.
    \]
\end{definition}
\begin{definition}
    At a~parameter $\theta$ where the parameter derivatives exist and belong to $C(K,\R^n)$, write
    \[ 
        D\mathcal R(\theta)[\dot\theta] 
        := 
        \left.\frac{d}{dt} \mathcal R(\theta+t\dot\theta)\right|_{t=0}
        =
        \sum_{i\in[P]} \dot\theta_i \frac{\partial}{\partial\theta_i}u_\theta. 
    \]
    Define
    \[
        \Theta_{\mathrm{fr}}
        :=
        \left\{
        \theta\in\Theta:
        \left\{
        \frac{\partial}{\partial\theta_j}u_\theta
        \right\}_{j\in[P]}
        \text{ is linearly independent in }C(K, \R^n)
        \right\}
    \]
    For $X=(x_1,\dots,x_m)\in K^m$ and $\theta\in\Theta$, define
    \begin{equation}\label{eq:F_X_and_DF_X}
        F_X:\nu\mapsto (u_\nu(x_i))_{i\in[m]},
        \qquad
        DF_X(\theta)
        =
        \left(
        \frac{\partial}{\partial\theta_j}u_\theta(x_i)
        \right)_{i\in[m], j\in[P]},
    \end{equation}
    where the Jacobian is defined whenever $F_X$ is differentiable. Note that each $\partial_{\theta_j}u_\theta(x_i)$ is an~$n$-dimensional column block, thus $DF_X(\theta)\in\R^{mn\times P}$.
\end{definition}

\subsection{Identification from finite samples}
We next introduce semianalytic and subanalytic sets, which provide the geometric framework for our general results. The use of~these notions is essential even for our main example: $\Theta_{\mathrm{reg}}^{\tanh}$ is semialgebraic, and hence semianalytic, because its defining conditions are finite combinations of~polynomial inequalities.
\begin{definition}[{\cite[Section 1.1]{valette2025subanalytic}}]\label{def:semi/sub-analyticity}
    A subset $E\subset\R^n$ is called
    \begin{itemize}
        \item \emph{semianalytic} if it is locally defined by~finitely many real analytic equalities and inequalities,
        \item \emph{globally semianalytic} if $V_n(E)$ is a~semianalytic subset of~$\R^n$, where $V_n\colon \R^n \rightarrow (-1,1)^n$ is the homeomorphism defined by
        \[
            V_n(\mathbf{z})
            =
            \frac{\mathbf{z}}{\sqrt{1+|\mathbf{z}|^2}}
            \qquad
            \text{for } \mathbf{z}\in\R^n,
        \]
        \item \emph{globally subanalytic} if it can be presented as a~projection of~a~globally semianalytic set onto the first $n$ coordinates,
        \item $\sigma$-\emph{subanalytic} if it is a~countable union of~globally subanalytic subsets of~$\R^n$.
    \end{itemize}
    Finally, let $E\subseteq \R^n$ be a~$\sigma$-subanalytic set. Then $f\colon E\rightarrow \R^m$ is a~$\sigma$-subanalytic function if its graph is a~$\sigma$-subanalytic set in $\R^{n+m}$. 
\end{definition}

\begin{figure}
    \centering
    \includegraphics[width=0.65\linewidth]{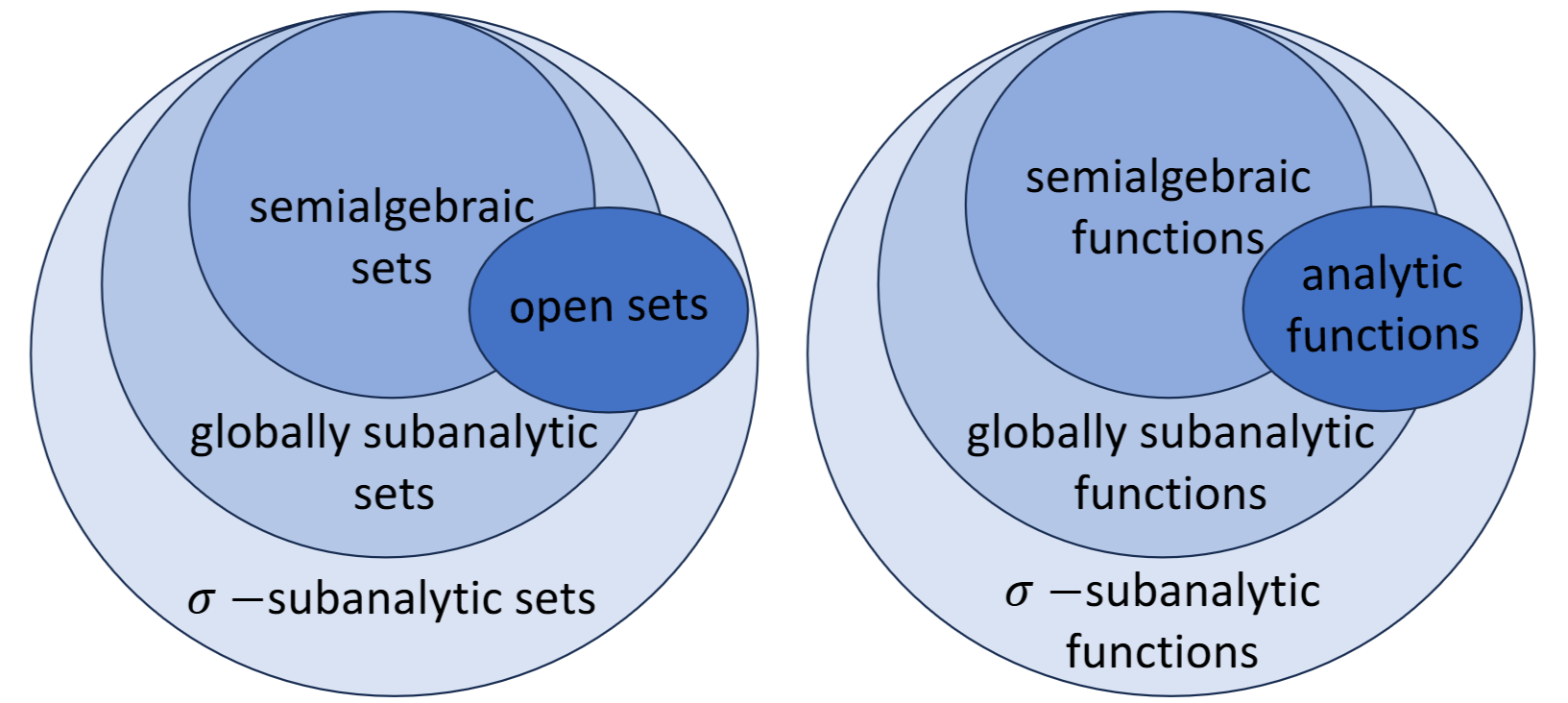}
    \caption{Inclusions between different classes of~sets; adapted from \cite[Figure~3]{valette2025subanalytic}.}
    \label{fig:subanalytic_semianalytic_analytic_sets}
\end{figure}

For a~more detailed discussion of~semianalytic, subanalytic, and $\sigma$-subanalytic sets and functions, see \cite{valette2025subanalytic} or \cite[Sections~1--3]{amir2023neural}. The inclusions between these classes of~sets are illustrated in Figure~\ref{fig:subanalytic_semianalytic_analytic_sets}. For this large variations of~sets, Amir et al. \cite{amir2023neural} estabilished a~general Finite witness theorem for identification a~mapping on the whoole domain only from finite number of~samples. This is a~corestone of~our general results and therefore we state it in the full statement below.

\begin{theorem}[{\cite[Theorem A.2]{amir2023neural}}]\label{thm:fwt}
    Let $\mathbb{M} \subseteq \R^p$, $\mathbb{W} \subseteq \R^q$ be $\sigma$-subanalytic sets of~dimension $D$ and $D_{\mathbf{\theta}}$ respectively. Let $F:\mathbb{M} \times \mathbb{W} \to \R$ be a~$\sigma$-subanalytic function. Define the set
    \[
        \mathcal{N}=\{\mathbf{z}\in \mathbb{M} \ | \ F(\mathbf{z};\mathbf{\theta})=0,\ \forall \mathbf{\theta}\in \mathbb{W} \}.
    \]
    Suppose that for all $\mathbf{z}\in \mathbb{M} \setminus \mathcal{N}$
    \begin{equation}\label{eq:assumption}
        \dim_H\{\mathbf{\theta}\in \mathbb{W} \ | \ F(\mathbf{z};\mathbf{\theta})= 0 \}\leq D_{\mathbf{\theta}}-1.\end{equation}
    Then for all $\left( \mathbf{\theta}^{(1)},\ldots, \mathbf{\theta}^{(D+1)}\right) \in \mathbb{W}^{D+1} \setminus \mathbb{F}$, where $\mathbb{F} \subset \mathbb{W}^{D+1}$ is a~set of~Hausdorff dimension strictly lower than $\dim_H(\mathbb{W}^{D+1})$, it holds that
    \begin{equation*}
        \mathcal{N}=\{ \mathbf{z}\in \mathbb{M} \ | \ F( \mathbf{z}; \mathbf{\theta}^{(i)})=0,\ \forall i=1,\ldots, D+1\} .
    \end{equation*}
    Moreover, if $\mathbb{W}$ is an~open and connected subset of~$\R^q$, and $F( \mathbf{z}; \mathbf{\theta})$ is analytic as a~function of~$ \mathbf{\theta}$ for all fixed $ \mathbf{z} \in \mathbb{M}$, then condition \eqref{eq:assumption} is not required, as it is automatically satisfied.
\end{theorem}

Thus, this result simply implies a~global identification of~two mappings on a~whole domain only from finite number of~samples. Although this result does not imply an~identifiability of~the mappings, it is a~key result for our finite samples identifiability result of~$\tanh$ networks.

\begin{lemma}\label{lem:generic_global_separation}
    Assume Setting~\ref{set:parmodel}. Suppose that
    \begin{enumerate}
        \item $\Theta$ is $\sigma$-subanalytic,
        \item $W\subset\mathbb R^d$ is a~nonempty connected open set,
        \item the mapping $u_\theta(x)$ from $\Theta \times W$ to $\R^n$ is $\sigma$-subanalytic,
        \item for every $\theta\in\Theta$, the function $x\mapsto u_\theta(x)$ is real analytic on $W$.
    \end{enumerate}
    Then, for Lebesgue almost every $(x_1,\dots,x_{2P+1})\in W^{2P+1}$,
    \begin{equation}\label{eq:generic_global_separation}
        u_\mu(x_i)=u_\nu(x_i)
        \text{ for all }i\in[2P+1]
        \quad\Longleftrightarrow\quad
        u_\mu=u_\nu\ \text{in }C(W, \R^n)
    \end{equation}
    holds simultaneously for all $\mu,\nu\in\Theta$.
\end{lemma}
\begin{proof}
    Set
    \[
        \mathbb M=\Theta \times \Theta,
        \qquad
        F((\mu, \nu);x)= \|u_\mu(x) - u_\nu(x)\|_2^2.
    \]
    For every $(\mu,\nu)\in\Theta \times \Theta$, the map $x\mapsto \|u_\mu(x)-u_\nu(x)\|_2^2$ is real analytic on $W$. Thus Theorem~\ref{thm:fwt} applies with $\mathbb W=W$. Since $D:=\dim_H(\mathbb M)\leq 2P$, adding further points preserves the equivalence. Hence, for almost every $(x_1,\dots,x_{2P+1})\in W^{2P+1}$,
    \[
        u_\mu(x_i)=u_\nu(x_i)\quad\text{for all }i\in[2P+1]
        \quad\Longleftrightarrow\quad
        u_\mu(x)=u_\nu(x)\quad\text{for all }x\in W
    \]
    simultaneously for all $(\mu,\nu)\in\Theta \times \Theta$.
\end{proof}

\subsection{Stable sampling}

The main general result of~this section is generalized Theorem \ref{thm:stable_sampling_tanh} for general activation function. As discussed above, the activation has to be analytic. The remaining task then remains to find a~set of~parameters, for which the network has full rank.

\begin{lemma}\label{lem:gamma_full_rank_cond}
Assume Setting \ref{set:parmodel}. For every $\gamma\in\Gamma$,
\[
    \theta\in\Theta_{\mathrm{fr}}
    \quad\Longleftrightarrow\quad
    \gamma(\theta)\in\Theta_{\mathrm{fr}}.
\]
In particular, for every compact
$\Xi\subset\Theta_{\mathrm{fr}}$, the set
\[
    \widetilde{\Xi}
    :=
    \bigcup_{\gamma\in\Gamma}\gamma(\Xi)
\]
is a~compact and $\Gamma$-invariant subset of
$\Theta_{\mathrm{fr}}$ containing $\Xi$.
\end{lemma}

\begin{proof}
Let $\theta\in\Theta_{\mathrm{fr}}$. Hence there exists a~point set $X=\{x_1,\dots,x_P\} \subset K$ such that
\[
\det\left(
\frac{\partial}{\partial\theta_j}u_\theta(x_i)
\right)_{i,j\in[P]}
\neq0.
\]
Recall \eqref{eq:F_X_and_DF_X}. Since $\mathcal R$ is $\Gamma$-invariant by~\eqref{eq:gamma_invariant},
\[
F_X\circ\gamma=F_X.
\]
Differentiating gives
\[
DF_X(\gamma(\theta))\,D\gamma(\theta)=DF_X(\theta).
\]
Because $\gamma$ is a~$C^1$-diffeomorphism, $D\gamma(\theta)$ is invertible, and therefore
\[
DF_X(\gamma(\theta))
=
DF_X(\theta)\bigl(D\gamma(\theta)\bigr)^{-1}.
\]
Hence
\[
\det DF_X(\gamma(\theta))
=
\det DF_X(\theta)\,
\det\bigl(D\gamma(\theta)\bigr)^{-1}
\neq0.
\]
Thus $\gamma(\theta)\in\Theta_{\mathrm{fr}}$. The converse implication follows by~applying the same argument to the inverse diffeomorphism $\gamma^{-1}\in\Gamma$.

Now let $\Xi\subset\Theta_{\mathrm{fr}}$ be compact and define
\[
    \widetilde{\Xi}
    :=
    \bigcup_{\gamma\in\Gamma}\gamma(\Xi).
\]
By~the first part, $\widetilde{\Xi}\subset\Theta_{\mathrm{fr}}$. Since every $\gamma\in\Gamma$ is continuous, each $\gamma(\Xi)$ is compact, and since $\Gamma$ is finite, $\widetilde{\Xi}$ is compact. Moreover, for every $\gamma_0\in\Gamma$,
\[
    \gamma_0(\widetilde{\Xi})
    =
    \bigcup_{\gamma_1\in\Gamma}
    (\gamma_0\circ\gamma_1)(\Xi)
    =
    \bigcup_{\gamma_2\in\Gamma}\gamma_2(\Xi)
    =
    \widetilde{\Xi},
\]
because $\{\gamma_0\circ\gamma_1:\gamma_1\in\Gamma\}=\Gamma$. Thus $\widetilde{\Xi}$ is $\Gamma$-invariant. Finally, $\Xi\subset\widetilde{\Xi}$ since $\operatorname{Id}\in\Gamma$.
\end{proof}
\begin{theorem}\label{thm:generic_local_sampling}
    Assume Setting~\ref{set:parmodel}. Suppose that
    \begin{enumerate}
        \item $\Theta_{\mathrm{fr}}$ is $\sigma$-subanalytic,
        \item $Y\supset K$ is connected and open, and $W\subset \operatorname{int} K$ is nonempty and open,
        \item $D\mathcal R(\theta)[v](x)$ is $\sigma$-subanalytic on $\Theta_{\mathrm{fr}}\times\R^P\times Y$,
        \item for every $(\theta,v)\in\Theta_{\mathrm{fr}}\times\R^P$, the function
    $x\mapsto D\mathcal R(\theta)[v](x)$ is real analytic on $Y$,
        \item for every $x\in K$, the map $\nu\mapsto u_\nu(x)$ is $C^1$ in a~neighborhood of~each $\theta\in\Theta_{\mathrm{fr}}$.
    \end{enumerate}
    Then, for Lebesgue almost every
    $(x_1,\dots,x_{2P+1})\in W^{2P+1}$, it holds
    that for every $\theta\in\Theta_{\mathrm{fr}}$, there exist $r_\theta,\delta_\theta>0$ such that
    \begin{equation}\label{eq:stablesampling_2P+1}
        \forall \mu,\nu\in B(\theta,r_\theta):\ \max_{i\in [2P+1]}\|u_\mu(x_i) - u_\nu(x_i)\|_2\geq \delta_\theta\|u_\mu - u_\nu\|_{L^\infty(K)}.
    \end{equation}
\end{theorem}
\begin{proof}
    Set
    \[
        F((\theta, v);x) := \|D\mathcal R(\theta)[v](x)\|_2^2
    \]
    and
    \[
        \mathbb M:=\Theta_{\mathrm{fr}}\times\R^P,
        \qquad
        \mathcal N
        :=
        \left\{
        (\theta,v)\in\mathbb M:
        F((\theta, v);x) =0
        \text{ for every }x\in Y
        \right\}.
    \]
    If $(\theta,v)\in\mathcal N$, then $x\mapsto F((\theta,v);x)$ vanishes on $K\subset Y$, so $v=0$ by~the definition of~$\Theta_{\mathrm{fr}}$. Conversely, if $v=0$, then the real analytic function $x\mapsto F((\theta,0);x)$ vanishes on $K$. Since $K$ contains the nonempty open set $W$ and $Y$ is connected, the identity theorem gives $(\theta,0)\in\mathcal N$. Thus
    \[
        \mathcal N=\Theta_{\mathrm{fr}}\times\{0\},
    \]
    and we apply Theorem~\ref{thm:fwt} with $\mathbb W=Y$. Since $D:=\dim_H(\mathbb{M})\leq 2P$, adding further points preserves theorem statement. Restricting to $W$, we have that for almost every $(x_1,\dots,x_{2P+1})\in W^{2P+1}$,
    \begin{equation}\label{eq:full_rank_from_finite_m}
        F((\theta, v);x_i)=0
        \text{ for all }i\in[2P+1]
        \quad\iff\quad
        v=0
    \end{equation}
    simultaneously for all $(\theta,v)\in\mathbb M$.

    Fix such a~tuple $X$. Then \eqref{eq:full_rank_from_finite_m} implies
    \begin{equation}\label{eq:full_rank_DF}
        \ker DF_X(\theta)=\{0\},
        \qquad \theta\in\Theta_{\mathrm{fr}},
    \end{equation}
    and hence $\operatorname{rank}DF_X(\theta)=P$.

    Now fix $\theta\in\Theta_{\mathrm{fr}}$. Since $\operatorname{rank}DF_X(\theta)=P$, there exists $I_\theta\subset[2P+1] \times [n]$, with $|I_\theta|=P$, such that the mapping
    \[
        F_{X_\theta}:\nu\longmapsto\left( \left(u_\nu(x_i)\right)_k \right)_{(i,k)\in I_\theta} \in \R^P
    \]
    satisfies
    \[
        \det DF_{X_\theta}(\theta)
        =
        \det\left(
            \frac{\partial}{\partial\theta_j} \left(u_\theta(x_i)\right)_k
        \right)_{(i,k)\in I_\theta,\,j\in[P]}
        \neq 0.
    \]
    By~the assumed $C^1$ regularity and the inverse function theorem, there exist open neighborhoods $U_\theta$ of~$\theta$ and $V_\theta$ of~$F_{X_\theta}(\theta)$ such that $F_{X_\theta}:U_\theta\to V_\theta$ is bijective with a~$C^1$ inverse.

    Choose $\rho_\theta>0$ such that $\overline{B(F_{X_\theta}(\theta),\rho_\theta)}
    \subset V_\theta$. Set $s_\theta:=\sigma_{\min}(DF_{X_\theta}(\theta))>0$. By the inverse function theorem,
    \[
        \bigl\|D(F_{X_\theta}^{-1})(F_{X_\theta}(\theta))\bigr\|_2
        =\bigl\|(DF_{X_\theta}(\theta))^{-1}\bigr\|_2
        =s_\theta^{-1}.
    \]
    By continuity, after decreasing $\rho_\theta$ if necessary,
    \[
        \sup_{y\in B(F_{X_\theta}(\theta),\rho_\theta)}
        \|D(F_{X_\theta}^{-1})(y)\|_2
        \le \frac{2}{s_\theta}.
    \]
    Since this ball is convex, $F_{X_\theta}^{-1}$ is Lipschitz there with constant $C_\theta:=2/s_\theta$. By continuity of~$F_{X_\theta}$, choose $r_\theta>0$ sufficiently small that
    \[
        B(\theta,r_\theta)\subset U_\theta,
        \qquad
        F_{X_\theta}(B(\theta,r_\theta))
        \subset B(F_{X_\theta}(\theta),\rho_\theta).
    \]
    Consequently, for all $\mu,\nu\in B(\theta,r_\theta)$,
    \[
        |\mu-\nu|
        \leq C_\theta\|F_{X_\theta}(\mu)-F_{X_\theta}(\nu)\|_2
        \leq C_\theta\sqrt P
        \max_{i\in[2P+1]}\|u_\mu(x_i)-u_\nu(x_i)\|_2.
    \]
    Finally, by~the local Lipschitz property \eqref{eq:lipschitz}, after possibly decreasing $r_\theta$, there exists $L_\theta>0$
    such that
    \[
        \|u_\mu-u_\nu\|_{L^\infty(K)}
        \leq L_\theta|\mu-\nu|
        \qquad
        \text{for all }\mu,\nu\in B(\theta,r_\theta).
    \]
    Combining these estimates yields
    \[
        \max_{i\in[2P+1]}\|u_\mu(x_i)-u_\nu(x_i)\|_2
        \geq
        \frac{1}{L_\theta C_\theta\sqrt P}
        \|u_\mu-u_\nu\|_{L^\infty(K)}.
    \]
    Taking
    \begin{equation}\label{eq:delta_theta}
        \delta_\theta=s_\theta/(2L_\theta\sqrt P)
    \end{equation}
    proves the claim.
\end{proof}
%
\begin{theorem}\label{thm:generic_global_sampling}
    Assume Setting~\ref{set:parmodel}. Suppose that
    \begin{enumerate}
        \item $\Theta_{\mathrm{fr}}$ is $\sigma$-subanalytic,
        \item $Y\supset K$ is connected and open, and $W\subset \operatorname{int} K$ is nonempty and open,
        \item $u_\theta(x)$ and $D\mathcal R(\theta)[v](x)$ are $\sigma$-subanalytic on $\Theta_{\mathrm{fr}}\times Y$ and $\Theta_{\mathrm{fr}}\times\R^P\times Y$, respectively,
        \item for every $(\theta,v)\in\Theta_{\mathrm{fr}}\times\R^P$, the functions $x\mapsto u_\theta(x)$ and $x\mapsto D\mathcal R(\theta)[v](x)$ are real analytic on $Y$.
        \item for every $x\in K$, the map $\nu\mapsto u_\nu(x)$ is $C^1$ in a~neighborhood of~each $\theta\in\Theta_{\mathrm{fr}}$.
    \end{enumerate}
    Suppose furthermore that global identifiability for $(\mathcal R,\Theta_{\mathrm{fr}})$ modulo $\Gamma$ holds.

    Then, for Lebesgue almost every $X=(x_1,\dots,x_{2P+1})\in W^{2P+1}$, it holds that for every compact subset $\Xi\subset\Theta_{\mathrm{fr}}$, there exists $\delta_\Xi>0$ such that
    \begin{equation}\label{eq:stablesamplingglobal_2P+1}
        \mbox{for all }\mu,\nu\in \Xi:\ \max_{i\in [2P+1]}\|u_\mu(x_i) - u_\nu(x_i)\|_2\geq \delta_\Xi \|u_\mu - u_\nu\|_{L^\infty(K)}.
    \end{equation}
\end{theorem}
\begin{proof}
   By~Theorem~\ref{thm:generic_local_sampling} and Lemma~\ref{lem:generic_global_separation}, the latter applied to the restricted family $\mathcal R|_{\Theta_{\mathrm{fr}}}$, there exist full-measure subsets $G_{\mathrm{loc}},G_{\mathrm{sep}} \subset W^{2P+1}$ on which the respective conclusions hold. Since the intersection of~finitely many full-measure sets has full measure, $G:=G_{\mathrm{loc}}\cap G_{\mathrm{sep}}$ has full measure in $W^{2P+1}$. Fix such a~tuple $X=(x_1,\dots,x_{2P+1})\in G$ and define
    \[
        q_X(\mu,\nu)
        :=
        \max_{i\in[2P+1]}
        \|u_\mu(x_i)-u_\nu(x_i)\|_2.
    \]
    Fix a~compact subset $\Xi\subset\Theta_{\mathrm{fr}}$. By~Lemma~\ref{lem:gamma_full_rank_cond}, there exists a~compact $\Gamma$-invariant set $\widetilde{\Xi}\subset\Theta_{\mathrm{fr}}$ with $\Xi\subset\widetilde{\Xi}$. Hence, in the remainder of~the proof, assume that $\Xi$ is $\Gamma$-invariant.
    
    By~Theorem~\ref{thm:generic_local_sampling}, for every $\theta\in\Xi$ there exist $r_\theta,\delta_\theta>0$ such that
    \begin{equation}\label{eq:generic_local_patch}
        q_X(\mu,\nu)
        \geq
        \delta_\theta
        \|u_\mu-u_\nu\|_{L^\infty(K)}
    \end{equation}
    whenever $\mu,\nu\in B(\theta,r_\theta)$.

    Since $\Xi$ is compact, there exists a~finite subset $A\subset\Xi$ with
    \begin{equation}\label{eq:Xi_cover}
        \Xi\subset
        \bigcup_{\theta\in A}B(\theta,r_\theta/2).
    \end{equation}
    Set
    \[
        r:=\frac12\min_{\theta\in A}r_\theta>0,
        \qquad
        \delta_1:=\min_{\theta\in A}\delta_\theta>0.
    \]
    For $\mu,\nu\in\Xi$ define
    \begin{equation}\label{eq:d_Gamma}
        d_\Gamma(\mu,\nu)
        :=
        \min_{\gamma_1,\gamma_2\in\Gamma}
        |\gamma_1(\mu)-\gamma_2(\nu)|.
    \end{equation}
    Since $\Gamma$ is finite, the infimum is attained. Assume that $d_\Gamma(\mu,\nu)<r$ and let $\gamma_1,\gamma_2\in\Gamma$ satisfy $d_\Gamma(\mu,\nu) = |\gamma_1(\mu)-\gamma_2(\nu)|$. By~the $\Gamma$-invariance of~$\Xi$ and the covering property \eqref{eq:Xi_cover}, there exists
    $\theta\in A$ such that
    \[
    \gamma_1(\mu)\in B(\theta,r_\theta/2).
    \]
    Moreover,
    \[
    |\gamma_2(\nu)-\theta|
    \leq
    |\gamma_2(\nu)-\gamma_1(\mu)|
    +
    |\gamma_1(\mu)-\theta|
    <
    r+\frac{r_\theta}{2} \le r_\theta.
    \]
    Therefore,
    \[
    \gamma_1(\mu),\gamma_2(\nu)\in B(\theta,r_\theta),
    \]
    and the local estimate \eqref{eq:generic_local_patch} applies. Using the invariance property \eqref{eq:gamma_invariant}, we conclude that
    \begin{equation}\label{eq:generic_near} 
    \mbox{for all }\mu,\nu\in \Xi :\ d_\Gamma(\mu,\nu)<r\Rightarrow q_X(\mu, \nu) \geq \delta_1\|u_\nu - u_\mu\|_{L^\infty(K)}.
    \end{equation}

    It remains to treat the case of~parameter pairs $\mu,\nu$ with $d_\Gamma(\mu,\nu)\geq r$. To this end, define
    \[
        B_r:=\{(\nu,\mu)\in \Xi\times \Xi:\ d_\Gamma(\nu,\mu)\geq r\}.
    \]
    We show that $B_r$ is compact. Since $\Gamma$ is finite,
    \[
    B_r=
    \bigcap_{\gamma_1,\gamma_2\in\Gamma}
    \{(\nu,\mu)\in\Xi\times\Xi:
    |\gamma_1(\nu)-\gamma_2(\mu)|\ge r\}.
    \]
    Since every $\gamma\in\Gamma$ is a~$C^1$ diffeomorphism, the map $(\nu,\mu)\mapsto |\gamma_1(\nu)-\gamma_2(\mu)|$ is continuous. Hence every set in the intersection is closed in $\Xi\times\Xi$. Since $\Xi\times\Xi$ is compact, these sets are compact, and therefore their finite intersection $B_r$ is compact.

    If $B_r=\emptyset$, then \eqref{eq:generic_near} proves the result with $\delta_\Xi=\delta_1$. Otherwise, the mapping
    \[
        (\mu,\nu)\longmapsto q_X(\mu,\nu)
    \]
    is continuous. Hence the minimum
    \[
        a
        :=
        \min_{(\mu,\nu)\in B_r}q_X(\mu,\nu)
    \]
    is attained. We claim that $a>0$. Indeed, if $a=0$, then for some $(\mu,\nu)\in B_r$ we have
    \[
        q_X(\mu,\nu)=0.
    \]
    By~Lemma~\ref{lem:generic_global_separation}, this implies $u_\mu=u_\nu$. Global identifiability modulo $\Gamma$ therefore yields $\mu\sim_\Gamma\nu$ on $Y$, and hence
    \[
        d_\Gamma(\mu,\nu)=0,
    \]
    contradicting $(\mu,\nu)\in B_r$. Thus $a>0$.

    Moreover, the mapping
    \[
        (\mu,\nu)\longmapsto
        \|u_\mu-u_\nu\|_{L^\infty(K)}
    \]
    is continuous. Hence the maximum
    \[
        M
        :=
        \max_{(\mu,\nu)\in B_r}
        \|u_\mu-u_\nu\|_{L^\infty(K)}
    \]
    is attained and finite. In fact $M>0$, since otherwise $u_\mu=u_\nu$ for every $(\mu,\nu)\in B_r$, contradicting $d_\Gamma(\mu,\nu)\geq r>0$ by~global identifiability.

    Consequently, for every $(\mu,\nu)\in B_r$,
    \[
        q_X(\mu,\nu)
        \geq a
        \geq
        \frac{a}{M}
        \|u_\mu-u_\nu\|_{L^\infty(K)}.
    \]
    Combining this estimate with \eqref{eq:generic_near} proves \eqref{eq:stablesamplingglobal_2P+1} with
    \[
        \delta_\Xi
        :=
        \min\left\{
        \delta_1,\frac{a}{M}
        \right\}>0.
    \]
    By local Lipschitz continuity and finite covering \eqref{eq:Xi_cover} of~the compact $\Xi$, $\mathcal R|_\Xi$ is globally Lipschitz. Choose a~Lipschitz constant $L_\Xi\ge\max_{\theta\in A}L_\theta$ and set $s:=\min_{\theta\in A}s_\theta$. Then
    \[
        M
        \le \max_{(\mu,\nu)\in B_r}L_\Xi|\mu-\nu|
        \le L_\Xi\operatorname{diam}(\Xi).
    \]
    Using \eqref{eq:delta_theta} therefore gives
    \[
        \delta_\Xi\ge
        \frac1{L_\Xi}\min\left\{
            \frac s{2\sqrt P},
            \frac{a}{\operatorname{diam}(\Xi)}
        \right\}.
    \]
\end{proof}
\begin{remark}
    We would like to mention that stable sampling can also be obtained by~a finite covering of~$\Xi$, without applying Theorem~\ref{thm:fwt}. More precisely, for every compact, $\Gamma^{\tanh}$-invariant set $\Xi\subset\Theta_{\mathrm{reg}}^{\tanh}$, there exist $m\in\N$ and $x_1,\dots,x_m\in K$ such that the stable sampling inequality holds. However, the advantage of~this approach is limited, since the resulting bound on $m$ may grow exponentially in $P$.
\end{remark}

\subsection{Local reconstruction}

The final result of~Section~\ref{sec:gen_results} concerns the local reconstruction of~analytic neural networks by~projected gradient descent and projected Gauss--Newton. Using standard arguments from numerical analysis, we establish local convergence for networks of~arbitrary architecture whose sampled Jacobian has full column rank at the target parameter.

\begin{corollary}\label{cor:local_pgd_c1}
Assume Setting~\ref{set:parmodel} and the following:
\begin{enumerate}
    \item $\Theta_{\mathrm{fr}}$ is $\sigma$-subanalytic,
    \item $W\subset\operatorname{int}K$ and $Y\supset K$ are connected open sets,
    \item $D\mathcal R(\theta)[v](x)$ is $\sigma$-subanalytic,
    \item for every $(\theta,v)\in\Theta_{\mathrm{fr}}\times\mathbb R^P$,
    the function $x\mapsto D\mathcal R(\theta)[v](x)$
    is real analytic on $Y$,
    \item for every $x\in K$, the map $\nu\mapsto u_\nu(x)$
    is $C^1$ in a~neighborhood of~every
    $\theta\in\Theta_{\mathrm{fr}}$.\label{it:ass_C1}
\end{enumerate}
For $\theta_\star,\theta\in\Theta$ and $X\in W^{2P+1}$, define
\[
    L_{\theta_\star}(\theta)
    :=
    \frac12\|F_X(\theta)-F_X(\theta_\star)\|_2^2.
\]

Then, for almost every $X\in W^{2P+1}$, the following holds. For every nonempty compact set $\Xi\subset\Theta_{\mathrm{fr}}$, there exist $r,\alpha_0,\mu,A,B>0$ such that, for every $\theta_\star\in\Xi$, every nonempty closed convex set $C\subset\mathbb R^P$, every $\vartheta\in C\cap\Gamma(\theta_\star)$, every $\alpha\in(0,\alpha_0]$, and every $\theta_1\in C\cap B(\vartheta,r)$,
\begin{itemize}
\item the \textbf{projected gradient descent} iteration
\[
\theta_{k+1}
=P_C\bigl(\theta_k-\alpha\nabla L_{\theta_\star}(\theta_k)\bigr)
\]
satisfies, for all $k\ge1$ and $q_\alpha:=1-\alpha\mu\in(0,1)$,
\begin{equation*}
\|\theta_k-\vartheta\|_2
\le 
q_\alpha^{k-1}\|\theta_1-\vartheta\|_2,
\qquad \text{and} \qquad
L_{\theta_\star}(\theta_k)
\le 
Aq_\alpha^{2(k-1)}
\|\theta_1-\vartheta\|_2^2.
\end{equation*}

Moreover, for every $k\ge1$,
\[
\|u_{\theta_k}-u_{\theta_\star}\|_{L^\infty(K)}
\le Bq_\alpha^{k-1}\|\theta_1-\vartheta\|_2 .
\]

\item the \textbf{projected Gauss-Newton} iteration
\[
\theta_{k+1}
=P_C\left(\theta_k- \big(J(\theta_k)^\top J(\theta_k)\big )^{-1} \nabla L_{\theta_\star}(\theta_k) \right)
\]
is well-defined converges Q-superlinearly to $\vartheta\in C\cap\Gamma(\theta_\star)$, meaining that
\[
\lim_{k\to\infty}
\frac{\|\theta_{k+1}-\vartheta\|_2}
{\|\theta_k-\vartheta\|_2}=0,
\]
unless convergence occurs in finitely many steps.
\end{itemize}
\end{corollary}

\begin{proof}
Fix $X$ such that $J(\theta):=DF_X(\theta)$ has full column rank for every $\theta\in\Theta_{\mathrm{fr}}$, which holds for almost every $X$ by~Theorem~\ref{thm:generic_local_sampling}. By~Lemma~\ref{lem:gamma_full_rank_cond}, we may assume, without loss of~generality, that $\Xi$ is $\Gamma$-invariant.

We first prove the PGD statement. Define
\[
s:=\min_{\xi\in\Xi}\sigma_{\min}(J(\xi))>0,
\qquad
M:=\max_{\xi\in\Xi}\|J(\xi)\|_2>0.
\]
By~assumption~\ref{it:ass_C1}, $F_X$ is $C^1$ on an~open neighborhood $U$ of~$\Xi$. Since $\Xi$ is compact, there exists $\eta>0$ such that
\[
V:=\{\theta\in\mathbb R^P:
\operatorname{dist}(\theta,\Xi)\le\eta\}\subset U.
\]
The Jacobian $J$ is uniformly continuous on the compact set $V$. Hence there exists $r\in(0,\eta)$ such that, for all $\xi\in\Xi$ and $\theta\in B(\xi,r)$,
\begin{equation}\label{eq:pgd_J_continuity}
\|J(\theta)-J(\xi)\|_2
\le\delta:=\frac{s^2}{6M}.
\end{equation}

Fix $\theta_\star\in\Xi$ and $\vartheta\in C\cap\Gamma(\theta_\star)$. By~$\Gamma$-invariance, $\vartheta\in\Xi$ and
\begin{equation}\label{eq:invariance_of_L}
    L_{\vartheta}=L_{\theta_\star},
    \qquad
    \nabla L_{\vartheta}=\nabla L_{\theta_\star}.
\end{equation}
Let $\theta\in B(\vartheta,r)$ and write
\[
h:=\theta-\vartheta,\qquad
e:=F_X(\theta)-F_X(\vartheta),\qquad
g:=\nabla L_{\theta_\star}(\theta)=J(\theta)^\top e.
\]
By~the fundamental theorem of~calculus,
\[
e=H_\theta h,\qquad
H_\theta:=\int_0^1J(\vartheta+th)\,dt.
\]
Since $\vartheta+th\in B(\vartheta,r)$ for $t\in[0,1]$, equation~\eqref{eq:pgd_J_continuity} implies
\[
\begin{aligned}
\|H_\theta-J(\vartheta)\|_2
&\le\int_0^1
\|J(\vartheta+th)-J(\vartheta)\|_2\,dt
\le\delta,\\
\|H_\theta\|_2
&\le M+\delta\le2M,\\
\|J(\theta)\|_2
&\le M+\delta\le2M,
\end{aligned}
\]
where we used $\delta\le M$, which follows from $s\le M$.

Moreover,
\begin{align*}
\|J(\theta)^\top H_\theta-J(\vartheta)^\top J(\vartheta)\|_2 &\le
\|(J(\theta)-J(\vartheta))^\top H_\theta\|_2
+\|J(\vartheta)^\top(H_\theta-J(\vartheta))\|_2\\
&\le \delta\|H_\theta\|_2+M\delta\\
&\le3M\delta=\frac{s^2}{2}.
\end{align*}
Consequently, since $g=J(\theta)^\top H_\theta h$,
\begin{align*}
\langle g,h\rangle
&=\langle J(\vartheta)^\top J(\vartheta)h,h\rangle +\langle
(J(\theta)^\top H_\theta-J(\vartheta)^\top J(\vartheta))h,h
\rangle\\
&\ge s^2\|h\|_2^2-\frac{s^2}{2}\|h\|_2^2 = \frac{s^2}{2}\|h\|_2^2.
\end{align*}
Furthermore,
\[
\|g\|_2
\le
\|J(\theta)\|_2 \|H_\theta\|_2 \|h\|_2
\le
4M^2\|h\|_2,
\qquad
\|e\|_2
\le
2M\|h\|_2.
\]

Choose
\[
\alpha_0:=\frac1{M^2},
\qquad
q_\alpha:=1-\frac{\alpha s^2}{2}\in(0,1).
\]
Since all eigenvalues of~$J(\vartheta)^\top J(\vartheta)$ belong to $[s^2,M^2]$, for $\alpha\le\alpha_0$ we have
\[
\|I-\alpha J(\vartheta)^\top J(\vartheta)\|_2
\le 1-\alpha s^2.
\]
Since $\vartheta\in C$, nonexpansiveness of~projection gives
\begin{align*}
\|P_C(\theta-\alpha g)-\vartheta\|_2
&\le
\|\theta-\alpha g-\vartheta\|_2
=\|(I-\alpha J(\theta)^\top H_\theta)h\|_2\\
&\le
\bigl(
\|I-\alpha J(\vartheta)^\top J(\vartheta)\|_2
+\alpha\|J(\theta)^\top H_\theta
-J(\vartheta)^\top J(\vartheta)\|_2
\bigr)\|h\|_2\\
&\le
\left(1-\alpha s^2+\frac{\alpha s^2}{2}\right)
\|h\|_2
=q_\alpha\|h\|_2.
\end{align*}

Thus, whenever $\theta\in C\cap B(\vartheta,r)$, we also have $P_C(\theta-\alpha g)\in C\cap B(\vartheta,r)$. By~induction, all iterates remain in this neighborhood and
\[
\|\theta_k-\vartheta\|_2
\le 
q_\alpha^{k-1}\|\theta_1-\vartheta\|_2
\longrightarrow0.
\]
Finally, using the bound on $\|e\|_2$ established above and \eqref{eq:invariance_of_L}, we obtain
\begin{align*}
L_{\theta_\star}(\theta_k)
&=
\frac12\|F_X(\theta_k)-F_X(\vartheta)\|_2^2\\
&\le
2M^2\|\theta_k-\vartheta\|_2^2 
\le
2M^2q_\alpha^{2(k-1)} \|\theta_1-\vartheta\|_2^2
\longrightarrow0.
\end{align*}
Taking $\mu:=\frac{s^2}{2}$ and $A=2M^2$ proves the first two claims. The last follows from \eqref{eq:lipschitz}, the $\Gamma$-invariance of~$\mathcal R$, and compactness of~$\Xi$, after possibly shrinking $r$. 

The PGN convergence conclusion follows from the full-column-rank property \eqref{eq:full_rank_DF} and \cite[Remark~5.2]{izmailov2026local}. This completes the proof.
\end{proof}
\section{Tanh results}\label{sec:tanh_results}

We now specialize the results of~Section~\ref{sec:gen_results} to $\tanh$ networks. Their realization maps are real analytic in the inputs and parameters. Hence, the main technical step is to establish linear independence of~the parameter derivatives. We prove this at every regular parameter, so $\Theta^{\tanh}_{\mathrm{reg}}=\Theta^{\tanh}_{\mathrm{reg,fr}}$.

\subsection{Private poles}

The key ingredient in proving injectivity of~$D\mathcal R^{\tanh}$ is the private-pole lemma below. Building on the pole analysis of~\textcite{vlavcic2021affine}, we prove a~stronger conclusion about how poles are shared among neurons. For one-dimensional inputs, every nonempty subset of~neurons in a~given hidden layer contains a~neuron with infinitely many simple poles at which all other neurons in the subset are holomorphic. These poles allow us to isolate the contribution of~a~single neuron in the linearized realization map. The lemma thus provides the main step in the full-rank proof. Before stating the lemma, we recall the definitions of~cluster points and cluster depth.

\begin{definition}[{\cite[Definition~23]{vlavcic2021affine}}]
    For $E\subset\C$, let $C^1(E)$ denote the set of~cluster points of~$E$, and set
    \[
        C^0(E):=E,
        \qquad
        C^k(E):=C^1(C^{k-1}(E)).
    \]
    The cluster depth of~$E$ is defined by
    \[
        L_C(E):=\inf\{k\in\N\colon C^k(E)=\emptyset\},
        \qquad
        L_C(\emptyset)=0.
    \]
\end{definition}
\begin{lemma}\label{lem:private_poles}
    Assume Setting~\ref{set:deep_tanh_nets} with $d=1$ and let $\theta\in\Theta_{\mathrm{reg}}^{\tanh}$. Let $\ell\in[L]$ and $\emptyset\neq J\subseteq[n_\ell]$. Then there exist $i_\star\in J$, $p\in\widehat{\C}$, and a~sequence $(z_n)_{n\in\N}$ with $z_n\to p$ such that
    \begin{enumerate}
        \item for all $n\in\N$, $z_n$ is a~simple pole of~$h_{i_\star}^\ell$;
        \item for all $j\in J\setminus\{i_\star\}$ and $n\in\N$, $h_j^\ell$ is holomorphic at $z_n$.
    \end{enumerate}
    For $\ell=1$, $p=\infty$. If $\ell\ge2$, then $p\in\C$ and, for every $1\le m<\ell$, $k\in[n_m]$, and $j\in[n_{m+1}]$, both $h_k^m$ and $q_j^{m+1}$ are meromorphic at $p$.
\end{lemma}

\begin{proof}
    We first consider $\ell\geq2$. For every $j\in[n_\ell]$, regard $q_j^\ell$ as the realization of~the scalar-output subnetwork with parameters
    \[
        \theta_j^\ell=(\widetilde W^r,\widetilde b^r)_{r\in[\ell]},
    \]
    constructed as follows. Set $I_\ell:=\{j\}$ and $I_0:=[d]$, and recursively define
    \[
        I_{r-1}
        :=
        \left\{
            i\in[n_{r-1}] :
            w_{ki}^r\neq0
            \text{ for some }k\in I_r
        \right\},
        \qquad r=\ell,\ldots,2.
    \]
    Then set
    \[
        \widetilde W^r:=W^r_{I_r,I_{r-1}},
        \qquad
        \widetilde b^r:=b^r_{I_r},
        \qquad r\in[\ell].
    \]
    Note that we neither change the architecture of~the original network nor the realization of~$q_j^\ell$. We only consider the subnetwork obtained By~retaining those vertices and edges that influence $q_j^\ell$, as in \cite[Definition~20]{vlavcic2021affine}. One can then easily verify that $\theta_j^\ell$ is regular according to Definition~\ref{def:regular_parameters} and that the corresponding scalar-output network is a~$\tanh$-network in the sense of~Setting~\ref{set:deep_tanh_nets}, of~depth $\ell-1$ with realization $q_j^\ell$.
    
    Then let $D_j^\ell$ be the natural domain of~$q_j^\ell$; see \cite[Definition~11]{vlavcic2021affine}, let $P_j^\ell$ denote the set of~its simple poles, and set
    \begin{equation}\label{eq:E_j^l}
        E_j^\ell := (\C\setminus D_j^\ell)\setminus P_j^\ell,
    \end{equation}
    so that $E_j^\ell$ consists of~the essential singularities and poles of~higher order. By~\cite[Proposition~8(ii-iii)]{vlavcic2021affine},
    \[
        \C\setminus D_j^\ell=\overline{P_j^\ell},
        \qquad
        L_C\bigl(\overline{P_j^\ell}\bigr)=\ell-1.
    \]
    Moreover, \cite[Lemma~5(ii)]{vlavcic2021affine} yields
    \[
        \overline{P_j^\ell}
        =
        P_j^\ell\cup C^1(P_j^\ell),
    \]
    and hence
    \[
        E_j^\ell=C^1(P_j^\ell).
    \]
    Since $L_C(C^1(P_j^\ell)) = L_C(P_j^\ell)-1$ for $P_j^\ell\neq\emptyset$, \cite[Lemma~5(vi)]{vlavcic2021affine}  implies
    \[
        L_C(P_j^\ell)
        =
        L_C\bigl(\overline{P_j^\ell}\bigr)
        =
        \ell-1
        \qquad
        \text{and}
        \qquad 
        L_C(E_j^\ell)=\ell-2.
    \]
    
    Similarly, for $m<\ell$ and $k\in[n_m]$, let $F_k^m$ denote the set of~essential singularities and higher-order poles of~$h_k^m$. Regard $h_k^m$ as the realization of~the scalar-output subnetwork $\tilde q_1^{m+1}$ with output weight $1$ and output bias $0$. The network parameters can be defined as above, so they are again regular and the network has depth $m$. Denote by~$\tilde E_1^{m+1}$ and $\tilde P_1^{m+1}$ the sets as above, only for $\tilde q_1^{m+1}$. Applying \cite[Proposition~8]{vlavcic2021affine} and the same argument as above gives
    \[
        F_k^m=\tilde E_1^{m+1}=C^1(\tilde P_1^{m+1}),
        \qquad
        L_C(\tilde P_1^{m+1})=m,
        \qquad
        L_C(F_k^m)=m-1\le\ell-2.
    \]
    Define the bad set
    \[
        B:=
        \bigcup_{j=1}^{n_\ell}E_j^\ell
        \,\cup\,
        \bigcup_{\substack{1\leq m<\ell\\k\in[n_m]}}F_k^m .
    \]
    This is a~finite union of~sets of~clustering depth at most
    $\ell-2$. Hence, by~\cite[Lemma~5(vi)]{vlavcic2021affine},
    \begin{equation}\label{eq:bad_set_depth}
        L_C(B)\leq\ell-2.
    \end{equation}

    Fix any $i\in J$. Since
    $L_C(P_i^\ell)=\ell-1$, \eqref{eq:bad_set_depth} and monotonicity
    of~clustering depth imply that
    $P_i^\ell\not\subseteq B$. We may therefore choose
    \[
        p\in P_i^\ell\setminus B.
    \]
    This immediately proves the last assertion of~the lemma that for $\ell\ge2$, $h_k^m$ and $q_j^{m+1}$ are meromorphic at $p$ for all $m<\ell$, $k\in[n_m]$ and $j\in[n_{m+1}]$. Indeed, $p\notin F_k^m$ for every $m<\ell$, so each $h_k^m$ is either holomorphic or has a~simple pole at $p$, and is therefore meromorphic there. Moreover, by~\eqref{eq:deep_tanh_layers}, for every $m<\ell$ and $j\in[n_{m+1}]$, $q_j^{m+1}$ is an~affine combination of~$h_k^{m}$ for $k\in[n_m]$, and hence it is meromorphic at $p$ as well.

    Moreover, $p\notin E_j^\ell$ for every $j\in J$. Hence, by~\eqref{eq:E_j^l}, every $q_j^\ell$ is either holomorphic at $p$ or has a~simple pole there. Set
    \[
        I:=\{j\in J:q_j^\ell\text{ has a~pole at }p\}.
    \]
    Then $i\in I$, so $I\neq\emptyset$. If $j\in J\setminus I$, then $q_j^\ell$ is holomorphic at $p$. We claim that $h_j^\ell$ has no poles in a~sufficiently small punctured neighborhood of~$p$. Suppose, for contradiction, that there exists a~sequence of~poles $z_n\to p$ of~$h_j^\ell$. Then
    \[
    q_j^\ell(z_n)\in \pi\mathrm{i}\Bigl(\mathbb Z+\frac12\Bigr)
    \]
    for every $n$. Since $q_j^\ell$ is continuous at $p$, the sequence $q_j^\ell(z_n)$ converges to $q_j^\ell(p)$. As $\pi\mathrm{i}(\mathbb Z+\frac12)$ is discrete, $q_j^\ell(z_n)$ must therefore be eventually constant, say $q_j^\ell(z_n)=c$ for all sufficiently large $n$. Since the $z_n$ accumulate at $p$ and $q_j^\ell-c$ is holomorphic, the identity theorem implies that $q_j^\ell\equiv c$ in a~neighborhood of~$p$, and hence identically by~analytic continuation, contradicting the regularity of~the corresponding subnetwork.

    For $j\in I$, since $q_j^\ell$ has a~simple pole at $p$, we can write
    \begin{equation}\label{eq:expansion_of_q_near_p}
        q_j^\ell(z)
        =
        \frac{\beta_j}{z-p}
        +\gamma_j
        +\varepsilon_j(z-p),
    \end{equation}
    where $\beta_j\in\C\setminus\{0\}$, $\gamma_j\in\C$, and
    $\varepsilon_j$ is holomorphic near $0$ with $\varepsilon_j(0)=0$.

    By~Lemma~\ref{lem:simple_pole_local_inverse}, after restricting to a~sufficiently small neighborhood of~$p$, $q_j^\ell$ maps the corresponding punctured neighborhood one-to-one onto a~neighborhood of~infinity. Consequently, for every sufficiently large $|k_j|$, there is a~unique point $z_{j,k_j}$ near $p$ satisfying
    \begin{equation}\label{eq:poles_indexed_by_k}
        q_j^\ell(z_{j,k_j})
        =
        \pi\mathrm{i}\left(k_j+\frac12\right).
    \end{equation}
    These points satisfy $z_{j,k_j}\to p$ as $|k_j|\to\infty$ and are
    simple poles of~$h_j^\ell$.

    It remains to show that, for some $i_\star\in I$, infinitely many
    of~these poles are not shared by~the other $h_j^\ell$ for $j\in J, j\neq i_\star$. Thus, for $i,j\in I$, define
    \[
        A_{ij}
        :=
        \left\{
            k_i\in\Z:
            z_{i,k_i}=z_{j,k_j}
            \text{ for some }k_j\in\Z
        \right\},
    \]
    where the finitely many indices for which
    \eqref{eq:poles_indexed_by_k} has not yet been defined are ignored,
    and write
    \[
        \overline d(A)
        :=
        \limsup_{N\to\infty}
        \frac{\#(A\cap[-N,N])}{2N+1}
    \]
    for the upper density of~$A\subseteq\Z$. In particular, every finite set $A\subseteq\Z$ satisfies $\overline d(A)=0$.
    
    \medskip
    \begin{claim}\label{cl:w/beta=const}
        Let $i,j\in I$ satisfy
        \[
        \overline d(A_{ij})>0
        \qquad\text{or}\qquad
        \overline d(A_{ji})>0.
        \]
        Then,
        \[
            \operatorname{supp}\bigl(W^\ell_{i,:}\bigr)
            =
            \operatorname{supp}\bigl(W^\ell_{j,:}\bigr),
            \qquad
            \frac{W_{i,:}^\ell}{\beta_i}
            =
            \frac{W_{j,:}^\ell}{\beta_j}.
        \]
    \end{claim}

    \begin{proof}[Proof of~Claim \ref{cl:w/beta=const}]
    After interchanging $i$ and $j$ if necessary, we may assume that
    \[
        \overline d(A_{ij})>0.
    \]
    Set $\rho:=\frac{\beta_j}{\beta_i}$. Then, at every common pole $z_n=z_{i,k_i}=z_{j,k_j}$,
    \eqref{eq:expansion_of_q_near_p} gives
    \begin{equation}\label{eq:common_pole}
        \pi\mathrm{i}
        \left(
            \rho\left(k_i+\frac12\right)
            -
            \left(k_j+\frac12\right)
        \right)
        -
        (\rho\gamma_i-\gamma_j)
        =
        \rho\varepsilon_i(z_n-p)
        -
        \varepsilon_j(z_n-p).
    \end{equation}
    Dividing by~$k_i+\frac12$ along common poles with $|k_i|\to\infty$ shows that
    \[
        \rho
        -
        \frac{k_j+\frac12}{k_i+\frac12}
        = O\left(\frac{1}{|k_i|}\right).
    \]
    Hence $\rho\in\R$ and
    \begin{equation}\label{eq:k_asymptotic}
        k_j=\rho k_i+O(1).
    \end{equation}

    \begin{claim}\label{cl:uniformly_bounded_gaps_k_i}
        If $\delta:=\overline{d}(A_{ij})>0$, then there exist $C>0$ and infinitely many pairs of~distinct indices $k_i,k_i'\in A_{ij}$ for which $|k_i-k_i'|\le C$.
    \end{claim}
    \begin{proof}
        Set
        \[
        m_N:= \#(A_{ij}\cap[-N,N]).
        \]
        By~the definition of~$\limsup$,
        \[
        \forall \varepsilon>0, \forall N_0>0, \exists N\ge N_0\colon (\delta-\varepsilon)(2N+1) \le m_N.
        \]
        Set $\varepsilon:=\frac\delta2$. Let $a_1<\dots<a_{m_N}$ be the increasing enumeration of~$A_{ij}\cap[-N,N]$. Thus, for the sum of~the gaps, it holds that
        \[
        \sum_{l=1}^{m_N-1}a_{l+1}-a_l = a_{m_N}-a_1 \le 2N.
        \]
        Let's prove the claim by~contradiction. Then, for every $C>0$, there are only finitely many distinct pairs $k_i,k_i'\in A_{ij}$ with $|k_i-k_i'|\le C$. Fix $C:=8/\delta$. Since there are only finitely many such pairs, for all sufficiently large $N$, at most finitely many of~the consecutive gaps $a_{l+1}-a_l$ are at most $C$. Moreover, along the arbitrarily large $N$ satisfying
        \begin{equation}\label{eq:choice_of_m_N}
            \delta N +\frac{\delta}{2} \le m_N,
        \end{equation}
        we have $m_N\to\infty$. Hence, for all sufficiently large such $N$, at least half of~the gaps satisfy $a_{l+1}-a_l>C$. Thus,
        \begin{equation}\label{eq:condition_on_C}
            C\cdot \frac{m_N-1}{2}
            \le
            \sum_{l=1}^{m_N-1}a_{l+1}-a_l
            \le 
            2N.
        \end{equation}
        Hence $m_N\le \frac{4N}{C}+1$ and by~the choice of~$C$, we obtain
        \[
        m_N\le \frac{\delta N}{2}+1.
        \]
        For $N>2/\delta$, this is strictly smaller than left-hand side of~\eqref{eq:choice_of_m_N}, contradicting the choice of~$N$. Thus there exist $C>0$ and infinitely many distinct pairs $k_i,k_i'\in A_{ij}$ with $|k_i-k_i'|\le C.$
    \end{proof}

    By~Claim \ref{cl:uniformly_bounded_gaps_k_i}, there exist infinitely many pairs of~distinct indices $k_i,k_i'\in A_{ij}$ for which
    $|k_i-k_i'|$ is uniformly bounded. Hence, by~\eqref{eq:k_asymptotic}, the corresponding differences $|k_j-k_j'|$ are uniformly bounded as well. Therefore the integer pairs
    \[
        \bigl(k_i-k_i',\,k_j-k_j'\bigr)
    \]
    range over a~finite set. Consequently, there exist fixed integers $d\neq0$ and $e$ such that
    \[
        k_i-k_i'=d,
        \qquad
        k_j-k_j'=e
    \]
    for infinitely many pairs of~common poles.
    Subtracting \eqref{eq:common_pole} for these pairs gives
    \[
        \pi\mathrm{i}(\rho d-e)
        =
        \rho\bigl(
            \varepsilon_i(z_n-p)-\varepsilon_i(z_{n'}-p)
        \bigr)
        -
        \bigl(
            \varepsilon_j(z_n-p)-\varepsilon_j(z_{n'}-p)
        \bigr).
    \]
    Letting $z_n,z_{n'}\to p$ yields
    \[
        \rho=\frac{e}{d}\in\Q.
    \]

    Consequently, at every common pole,
    \[
        \rho\left(k_i+\frac12\right)
        -
        \left(k_j+\frac12\right)
        =
        \frac{2ek_i-2dk_j+e-d}{2d} \in \frac{1}{2d}\Z.
    \]
    Therefore the left-hand side of
    \eqref{eq:common_pole} belongs to the translated lattice
    \[
        S:=
        \frac{\pi\mathrm{i}}{2d}\Z
        -
        (\rho\gamma_i-\gamma_j).
    \]
    Any two distinct points of~$S$ have distance at least $\pi/(2d)$. On the other hand, by~\eqref{eq:common_pole}, these lattice points converge to $0$ along the common poles. Hence they must equal $0$ for all sufficiently large common poles: indeed, once their modulus is less than $\pi/(4d)$, the positive separation of~the lattice leaves at most one possible lattice point, and its limit is $0$. Thus
    \[
        \rho\varepsilon_i(z_n-p)
        =
        \varepsilon_j(z_n-p)
    \]
    for infinitely many $z_n\to p$. The identity theorem therefore gives
    \[
        \rho\varepsilon_i=\varepsilon_j,
        \qquad\text{or equivalently}\qquad
        \beta_i^{-1}\varepsilon_i
        =
        \beta_j^{-1}\varepsilon_j.
    \]
    Using \eqref{eq:expansion_of_q_near_p}, we therefore obtain
    \begin{equation}\label{eq:diff_q/beta_is_constant}
    \frac{q_i^\ell}{\beta_i}
    -
    \frac{q_j^\ell}{\beta_j}
    =
    \operatorname{const.}
    \end{equation}
    on a~punctured neighborhood of~$p$.
    
    Recall that $D_i^\ell$ and $D_j^\ell$ denote the natural domains of~$q_i^\ell$ and $q_j^\ell$, respectively, and set
    \[
    D:=D_i^\ell\cap D_j^\ell.
    \]    
    By~\cite[Proposition~8]{vlavcic2021affine}, both $\mathbb C\setminus D_i^\ell$ and $\mathbb C\setminus D_j^\ell$ are countable. Consequently, $\mathbb C\setminus D$ is countable, and hence $D$ is path-connected. Moreover, since the original network is well-defined
    for every real input, $\mathbb R\subset D$. The left-hand side of~\eqref{eq:diff_q/beta_is_constant} is holomorphic on $D$, so the identity theorem implies that \eqref{eq:diff_q/beta_is_constant} holds throughout $D$, and in particular on $\mathbb R$.
    
    Then \eqref{eq:diff_q/beta_is_constant} implies that on $\R$,
    \begin{equation}\label{eq:linear_dependence_parents}
        \sum_{s\in[n_{\ell-1}]}
        \left(
            \frac{W_{is}^\ell}{\beta_i}
            -
            \frac{W_{js}^\ell}{\beta_j}
        \right)
        h_s^{\ell-1}
        =
        C,
        \qquad C\in\R.
    \end{equation}
    Set
    \[
        c_s
        :=
        \frac{W_{is}^\ell}{\beta_i}
        -
        \frac{W_{js}^\ell}{\beta_j},
        \qquad s\in[n_{\ell-1}].
    \]
    We claim that every coefficient $c_s$ vanishes. Otherwise, since the
    functions $h_s^{\ell-1}$ are real-valued on $\R$, taking either the real
    or imaginary part of~\eqref{eq:linear_dependence_parents} yields a
    nontrivial real relation
    \begin{equation}\label{eq:subnetwork_c_s_vanishing}
        \sum_{s\in[n_{\ell-1}]} a_s h_s^{\ell-1}=A,
        \qquad a_s,A\in\R,
    \end{equation}
    with not all $a_s$ equal to zero. Let
    \[
        I':=\{s\in[n_{\ell-1}]:a_s\neq0\}.
    \]
    Consider the scalar-output subnetwork \eqref{eq:subnetwork_c_s_vanishing} with output weights $(a_s)_{s\in I'}$ and output bias $-A$. Similarly, one can verify that the parameters of~the network are regular and the network is non-trivial. However, its realization is identically zero, contradicting \cite[Proposition~1]{vlavcic2021affine}. Therefore
    \begin{equation}\label{eq:w/beta=const}
        \frac{W_{i,:}^\ell}{\beta_i}
        =
        \frac{W_{j,:}^\ell}{\beta_j}.
    \end{equation}
    In particular, combining it with \eqref{eq:every_neuron_has_parent}, we have
    \begin{equation}\label{eq:equality_of_supp_W}
        \operatorname{supp}(W_{i,:}^\ell)
        =
        \operatorname{supp}(W_{j,:}^\ell)
        \neq\emptyset.
    \end{equation}
    \end{proof}

    Now fix any $i_0\in I$ and set $G_0:=\{i_0\}$. For $m\geq0$, define
    \[
        G_{m+1}
        :=
        G_m\cup
        \left\{
            j\in I:
            \exists i\in G_m
            \text{ such that }
            \overline d(A_{ij})>0
            \text{ or }
            \overline d(A_{ji})>0
        \right\},
    \]
    and set
    \[
        G:=\bigcup_{m\geq0}G_m.
    \]
    By~Claim~\ref{cl:w/beta=const}, whenever $j\in G_{m+1}\setminus G_m$, there exists $i\in G_m$ such that \eqref{eq:w/beta=const} holds. Hence also \eqref{eq:equality_of_supp_W} holds, and choosing any $s\in\operatorname{supp}(W^\ell_{i,:})$ gives
    \begin{equation*}
        \frac{\beta_j}{\beta_i}
        =
        \frac{w_{js}^\ell}{w_{is}^\ell}
        \in\R\setminus\{0\}.
    \end{equation*}
    Iterating along the construction of~$G$ therefore yields
    \begin{equation}\label{eq:beta_j/beta_i_is_real}
        \frac{\beta_j}{\beta_i}\in\R\setminus\{0\},
        \qquad
        \frac{W_{i,:}^\ell}{\beta_i}
        =
        \frac{W_{j,:}^\ell}{\beta_j},
        \qquad \text{for every } i,j\in G.
    \end{equation}

    \medskip
    \begin{claim}\label{cl:pos_density_wrt_G}
    There exists $i_\star\in G$ such that, setting
    \[
        S_{i_\star}^G
        :=
        \left\{
        k\in\Z:
        z_{i_\star,k}
        \text{ is not a~pole of~}h_j^\ell
        \text{ for any }j\in G\setminus{i_\star}
        \right\},
    \]
    we have
    \[
        \overline d\left(S_{i_\star}^G
        \right)>0.
    \]
    \end{claim}

    \begin{proof}[Proof of~Claim \ref{cl:pos_density_wrt_G}]
    If $G=\{i_0\}$, set $i_\star=i_0$. Since $G\setminus\{i_\star\}=\emptyset$, every sufficiently large $|k|$ satisfies that $z_{i_\star,k}$ is not a~pole of~$h_j^\ell$ for any $j\in G\setminus\{i_\star\}$. Hence the set of~all such $k$ has upper density $1$ and the claim follows. Hence assume $|G|>1$. By~the preceding observation \eqref{eq:beta_j/beta_i_is_real}, there exist $\vartheta\in\R$ and $r_j\in\R\setminus\{0\}$ such that
    \[
        \beta_j=e^{\mathrm{i}\vartheta}r_j,
        \qquad j\in G.
    \]
    Moreover, for $s\in[n_{\ell-1}]$, we may define
    \[
        c_s:=\frac{w_{js}^\ell}{\beta_j},
    \]
    which is independent of~$j\in G$. Therefore
    \[
        q_j^\ell(z)
        =
        b_j^\ell
        +
        \beta_j\sum_{s\in[n_{\ell-1}]}c_s h_s^{\ell-1}(z)
        =
        r_jT(z)+b_j^\ell,
    \]
    where
    \[
        T(z)
        :=
        e^{\mathrm{i}\vartheta}
        \sum_{s\in[n_{\ell-1}]}c_s h_s^{\ell-1}(z).
    \]
    Since $q_j^\ell$ has a~simple pole at $p$ and $r_j\neq0$, $T$ has a~simple pole at $p$ as well. Hence, Lemma~\ref{lem:simple_pole_local_inverse} shows that, after restricting to a~sufficiently small neighborhood of~$p$, $T$ is one-to-one on the corresponding punctured neighborhood and maps it onto a~neighborhood of~infinity.

    In the variable $t=T(z)$, the poles of~$h_j^\ell$ are therefore,
    for sufficiently large modulus, in one-to-one correspondence with
    the points
    \begin{equation}\label{eq:pole_lattice_T}
        t
        =
        -\frac{b_j^\ell}{r_j}
        +
        \frac{\pi\mathrm{i}}{r_j}
        \left(k_j+\frac12\right),
        \qquad k_j\in\Z.
    \end{equation}
    Thus the poles of~each $h_j^\ell$ lie on the vertical line
    \[
        \Re t=-\frac{b_j^\ell}{r_j}.
    \]
    In particular, if the poles of~$h_i^\ell$ and $h_j^\ell$ lie on different vertical lines, then $h_i^\ell$ and $h_j^\ell$ cannot share a~pole.

    Thus, let
    \[
        \mathcal A_G
        :=
        \left\{
            \frac{b_j^\ell}{r_j}:j\in G
        \right\}.
    \]
    Choose any $\alpha\in\mathcal A_G$ and denote by
    \[
        G_\alpha
        :=
        \left\{
            j\in G:
            \frac{b_j^\ell}{r_j}=\alpha
        \right\}
    \]
    the set of~indices $j\in G$ for which the poles of~$h_j^\ell$ lie on the vertical line
    \[
        \Re t=-\alpha.
    \]
    Since $G_\alpha$ is finite and nonempty, we may choose
    \[
        i_\star\in G_\alpha
        \qquad\text{such that}\qquad
        |r_{i_\star}|
        =
        \max_{j\in G_\alpha}|r_j|.
    \]

    This maximizer is unique. Indeed, suppose that
    $j\in G_\alpha\setminus\{i_\star\}$ also satisfies
    \[
        |r_j|=|r_{i_\star}|.
    \]
    Then $r_j=\pm r_{i_\star}$. Since $i_\star,j\in G_\alpha$,
    \[
        b_j^\ell=\pm b_{i_\star}^\ell
    \]
    with the same sign. Moreover, \eqref{eq:w/beta=const} gives
    \[
        W_{j,:}^\ell=\pm W_{i_\star,:}^\ell
    \]
    with the same sign. Therefore
    \[
        (W_{j,:}^\ell,b_j^\ell)
        =
        \pm(W_{i_\star,:}^\ell,b_{i_\star}^\ell),
    \]
    contradicting \eqref{eq:regularity_of_params_2}.
    
    Now set
    \[
        H
        :=
        \left\{
            j\in G_\alpha\setminus\{i_\star\}:
            h_j^\ell\text{ and }h_{i_\star}^\ell
            \text{ share a~pole}
        \right\}.
    \]
    For $j\in H$, equality of~a~common pole in \eqref{eq:pole_lattice_T} gives
    \[
        \frac{k_{i_\star}+\frac12}{r_{i_\star}}
        =
        \frac{k_j+\frac12}{r_j}
    \]
    for some $k_{i_\star},k_j\in\Z$, and hence
    \[
        \frac{r_j}{r_{i_\star}}
        =
        \frac{2k_j+1}{2k_{i_\star}+1}
        \in\Q.
    \]
    Write
    \[
        \frac{r_j}{r_{i_\star}}
        =
        \frac{a_j}{d_j},
        \qquad
        \gcd(|a_j|,d_j)=1,
        \qquad d_j>0.
    \]
    By~uniqueness of~the maximizer, $|r_j|<|r_{i_\star}|$, and therefore
    \[
        |a_j|<d_j,
        \qquad d_j\ge2.
    \]
    Moreover, every common pole satisfies
    \[
        a_j(2k_{i_\star}+1)
        =
        d_j(2k_j+1).
    \]
    Since $a_j$ and $d_j$ are coprime, this implies the necessary condition
    \begin{equation}\label{eq:d_divides_odd}
        d_j\mid(2k_{i_\star}+1).
    \end{equation}
    
    It remains to choose a~positive-density set of~indices for which \eqref{eq:d_divides_odd} fails for every $j\in H$. If $d_j$ is even, then \eqref{eq:d_divides_odd} is impossible since $2k_{i_\star}+1$ is odd. For the remaining indices, set
    \[
        D
        :=
        \prod_{\substack{j\in H\\ d_j\ {\rm odd}}}d_j,
    \]
    with the convention $D=1$ if the product is empty. If $k_{i_\star}\in D\Z$, then every odd $d_j$ occurring in the product divides $k_{i_\star}$, and hence
    \[
        2k_{i_\star}+1\equiv1\pmod{d_j}.
    \]
    Thus \eqref{eq:d_divides_odd} fails for every $j\in H$.
    
    Consequently, for every sufficiently large $k_{i_\star}\in D\Z$, the point $z_{i_\star,k_{i_\star}}$
    \begin{itemize}
        \item is not a~pole of~$h_j^\ell$ for any $j\in H$;
        \item is not a~pole of~$h_j^\ell$ for any $j\in G_\alpha\setminus(H\cup\{i_\star\})$, by~the definition of~$H$;
        \item hence, is not a~pole of~$h_j^\ell$ for any $j\in G_\alpha\setminus\{i_\star\}$.
    \end{itemize}
    For $j\in G\setminus G_\alpha$, the sufficiently large poles of~$h_j^\ell$ lie on a~different vertical line \eqref{eq:pole_lattice_T} and therefore cannot coincide with $z_{i_\star,k_{i_\star}}$ for all sufficiently large $|k_{i_\star}|$. Since $D\Z$ has density $1/D$, we obtain
    \[
        \overline d\left(
            \left\{
                k\in\Z:
                z_{i_\star,k}
                \text{ is not a~pole of~}h_j^\ell
                \text{ for any }j\in G\setminus\{i_\star\}
            \right\}
        \right)
        \geq\frac1D>0.
    \]
    \end{proof}

    \begin{claim}\label{cl:pos_density_wrt_I}
        Let $i_\star\in G$ be as in Claim~\ref{cl:pos_density_wrt_G}. Then
        \[
            \overline d\left(
            \left\{
                k\in\Z:
                z_{i_\star,k}
                \text{ is not a~pole of~}h_j^\ell
                \text{ for any }j\in I\setminus\{i_\star\}
            \right\}
            \right)>0.
        \]
    \end{claim}
    
    \begin{proof}[Proof of~Claim \ref{cl:pos_density_wrt_I}]
        For every $j\in I\setminus G$,
        \[
            \overline d(A_{i_\star j})=0.
        \]
        Indeed, since $i_\star\in G$, we have $i_\star\in G_m$ for some $m\geq0$. If $\overline d(A_{i_\star j})>0$, then the definition of~$G_{m+1}$ would imply $j\in G$, a~contradiction. Since $I\setminus G$ is finite, setting $U := \bigcup_{j\in I\setminus G} A_{i_\star j}$ gives $\overline d(U)=0.$ Since
        \[
            \overline d(S_{i_\star}^G\setminus U)
            \ge
            \overline d(S_{i_\star}^G)-\overline d(U)
            =
            \overline d(S_{i_\star}^G)>0,
        \]
        removing $U$ from $S_{i_\star}^G$ preserves positive upper density. Since $S_{i_\star}^G\setminus U$ consists, up to finitely many indices, of~those $k$ for which $z_{i_\star,k}$ is not a~pole of~$h_j^\ell$ for any $j\in I\setminus\{i_\star\}$, we obtain
        \[
            \overline d\left(
            \left\{
            k\in\Z:
            z_{i_\star,k}
            \text{ is not a~pole of~}h_j^\ell
            \text{ for any }j\in I\setminus\{i_\star\}
            \right\}
            \right)>0.
        \]
    \end{proof}

    Since $|k_n|\to\infty$, we have $z_n\to p$. For $j\in J\setminus I$, the function $q_j^\ell$ is holomorphic at $p$, so $h_j^\ell$ has no poles in a~sufficiently small punctured
    neighborhood of~$p$. After discarding finitely many $z_n$, we therefore have
    \[
        h_j^\ell
        \quad\text{holomorphic at }z_n,
        \qquad
        j\in J\setminus\{i_\star\}.
    \]
    This proves the assertion for $\ell\geq2$.

    It remains to consider $\ell=1$. Since $d=1$, condition~\eqref{eq:every_neuron_has_parent} yields
    \[
        W^1_{j,:}=(w_j^1)\neq0,
        \qquad j\in[n_1],
    \]
    and hence $w_j^1\neq0$ for every $j\in[n_1]$. In this case
    \[
        q_j^1(z)=w_j^1z+b_j^1,
    \]
    so the poles of~$h_j^1$ are
    \[
        -\frac{b_j^1}{w_j^1}
        +
        \frac{\pi\mathrm{i}}{w_j^1}
        \left(k_j+\frac12\right),
        \qquad k_j\in\Z.
    \]
    Thus we are exactly in the situation of~\eqref{eq:pole_lattice_T}, with $t=z$, $r_j=w_j^1$, and $b_j^\ell=b_j^1$. Repeating the argument following \eqref{eq:pole_lattice_T}, choose a~vertical line and $h_{i_\star}^1$ with poles on it with maximal $|w_{i_\star}^1|$. Regularity of~parameters makes this maximizer unique, and the same divisibility argument produces infinitely many poles of~$h_{i_\star}^1$ which are not shared by~any other $h_j^1$ for $j\in J$. Since their indices satisfy $|k|\to\infty$, their moduli tend to infinity. This completes the proof.
\end{proof}

\begin{lemma}\label{lem:simple_pole_local_inverse}
    Let $p\in\C$ and let $f$ be meromorphic in a~neighborhood of~$p$. Suppose that $p$ is a~simple pole of~$f$. Then there exists a~neighborhood $U\subset\C$ of~$p$ such that $f$ is one-to-one on $U\setminus\{p\}$ and maps $U\setminus\{p\}$ onto a~neighborhood of~infinity.
\end{lemma}
\begin{proof}
    Since $p$ is a~simple pole of~$f$, there exist $\beta\in\C\setminus\{0\}$ and a~function $\varepsilon$ holomorphic near $p$ such that
    \[
        f(z)=\frac{\beta}{z-p}+\varepsilon(z).
    \]
    Hence, the function
    \[
        g(z):=\frac{1}{f(z)}
        =
        \frac{z-p}{\beta+(z-p)\varepsilon(z)}
    \]
    extends holomorphically to $p$ (because $\beta\neq0)$ and has a~simple zero there. Hence
    \[
        g(p)=0,
        \qquad
        g'(p)=\frac{1}{\beta}\neq0.
    \]
    By~the holomorphic inverse function theorem, after restricting to a~sufficiently small neighborhood $U$ of~$p$, $g$ is one-to-one and $g(U)$ is a~neighborhood of~$0$. Since $f=1/g$, the claim follows.
\end{proof}

\subsection{Full rank}

Using Lemma~\ref{lem:private_poles}, we now prove that the differential of the $\tanh$ realization map has full rank at every regular parameter $\theta$. Equivalently, the functions $\partial_{\theta_i} u_\theta$, $i\in[P]$, are linearly independent in $C^\infty(\R^d,\R^n)$. The lemma below treats the case of~one-dimensional inputs, $u_\theta\colon\R\to\R^n$. We then use this case to prove Theorem~\ref{thm:full_rank_realization_d} for $u_\theta\colon\R^d\to\R^n$.

\begin{lemma}\label{lem:full_rank_realization}
    Assume Setting~\ref{set:deep_tanh_nets} with $d=1$ and let $\theta\in\Theta_{\mathrm{reg}}^{\tanh}$. Then
    \[
        D\mathcal R(\theta)[\dot\theta]\equiv0 \quad \text{on } \R
        \qquad\Longrightarrow\qquad
        \dot\theta=0.
    \]
\end{lemma}

\begin{proof}
    Differentiating $h_i^\ell$ in the parameter direction $\dot\theta$ gives
    \begin{equation}\label{eq:variation_h_simple}
        \dot h_j^\ell
        =
        \bigl(1-(h_j^\ell)^2\bigr)\dot q_j^\ell .
    \end{equation}
    We use the same notation for the analytic continuations of~these functions to the complex plane.

    We first record the peeling argument that will be used at every layer. Suppose that, for some nonempty set $J\subset[n_\ell]$,
    \begin{equation}\label{eq:peeling_relation}
        0
        =
        c+
        \sum_{j\in J}
        \left(
            \dot a_j h_j^\ell+a_j\dot h_j^\ell
        \right),
        \qquad a_j\in\R.
    \end{equation}
    Apply Lemma~\ref{lem:private_poles} to $J$ and let $i_\star\in J$ be the resulting index satisfying that $h_{i_\star}^\ell$ has private simple poles $(z_n)_{n\in\N}$.
    
    If $\ell\geq2$, then $z_n\to p$ for some $p\in\C$, and all $h_k^m$ and $q_j^{m+1}$ are meromorphic at $p$ for $m<\ell$, $k\in[n_m]$, and $j\in[n_{m+1}]$. Since
    \[
        \dot q_i^1(z)
        =
        \dot b_i^1+\dot w_i^1 z,
    \]
    $\dot q_i^1$ is holomorphic on $\C$, and \eqref{eq:variation_h_simple} implies that $\dot h_i^1$ is meromorphic at $p$. Moreover, for $2\leq m<\ell$,
    \begin{equation}\label{eq:variation_q}
        \dot q_i^m
        =
        \dot b_i^m+
        \sum_{j\in[n_{m-1}]}
        \left(
            \dot w_{ij}^m h_j^{m-1}
            +
            w_{ij}^m\dot h_j^{m-1}
        \right).
    \end{equation}
    Hence, by~induction over $m<\ell$ and \eqref{eq:variation_h_simple}, $\dot q_i^m$ and $\dot h_i^m$ are meromorphic at $p$ for all $m<\ell$ and $i\in[n_m]$. Moreover, the same formula and argument also shows that $\dot q_j^\ell$ is meromorphic at $p$ for every $j\in[n_\ell]$. Hence, there exists a~neighborhood $U$ of~$p$ in which they are holomorphic except $p$. As $z_n\to p$ and $z_n\neq p$, after discarding finitely many terms we have $z_n\in U \setminus{p}$, and hence all these functions are holomorphic at $z_n$.
    
    If $\ell=1$, then $|z_n|\to\infty$, while $\dot q_j^1$ is holomorphic for every $j\in J$.
    
    In either case, for every $j\neq i_\star$, privacy of~$z_n$ implies that $h_j^\ell$ is holomorphic at $z_n$. Since $\dot q_j^\ell$ is holomorphic there as well, \eqref{eq:variation_h_simple} shows that
    $\dot h_j^\ell$ is holomorphic at $z_n$.
    
    Suppose first that $a_{i_\star}=0$. Then the only possible singular term in \eqref{eq:peeling_relation} at $z_n$ is
    \[
        \dot a_{i_\star}h_{i_\star}^\ell.
    \]
    Since $h_{i_\star}^\ell$ has a~pole at $z_n$, it follows that
    \[
        \dot a_{i_\star}=0.
    \]
    
    Suppose now that $a_{i_\star}\neq0$. Since $z_n$ is a~simple pole of~$h_{i_\star}^\ell$, the term $1-(h_{i_\star}^\ell)^2$ has a~double pole there. Therefore, if $\dot q_{i_\star}^\ell(z_n)\neq0$, then \eqref{eq:variation_h_simple} implies that $a_{i_\star}\dot h_{i_\star}^\ell$ has a~double pole at $z_n$. On the other hand, $\dot a_{i_\star}h_{i_\star}^\ell$ has at most a~simple pole, while all terms corresponding to $j\neq i_\star$ are holomorphic there. Hence the double pole cannot cancel, and
    \[
        \dot q_{i_\star}^\ell(z_n)=0
        \qquad
        \text{for all sufficiently large }n.
    \]

    If $\ell\geq2$, then $\dot q_{i_\star}^\ell$ is meromorphic at $p$ and has infinitely many zeros $z_n\to p$. Hence, for a~sufficiently large $k\geq0$, the function
    \[
        g(z):=(z-p)^k\dot q_{i_\star}^\ell(z)
    \]
    extends holomorphically to $p$. Since
    \[
        g(z_n)=(z_n-p)^k\dot q_{i_\star}^\ell(z_n)=0
    \]
    for infinitely many $z_n\to p$, the identity theorem gives $g\equiv0$. As $(z-p)^k\neq0$ for $z\neq p$, it follows that
    \[
        \dot q_{i_\star}^\ell\equiv0.
    \]
    For $\ell=1$, $\dot q_{i_\star}^1(z) = \dot w_{i_\star}^1z+\dot b_{i_\star}^1$ vanishes at infinitely many distinct private poles. Since a~nonzero affine function has at most one zero, it follows that
    \[
        \dot q_{i_\star}^\ell\equiv0,
        \qquad
        \dot h_{i_\star}^\ell\equiv0.
    \]

    Returning to \eqref{eq:peeling_relation}, the only remaining singular contribution of~the $i_\star$-th summand is $\dot a_{i_\star}h_{i_\star}^\ell$. Since $h_{i_\star}^\ell$ has a~pole, necessarily
    \[
        \dot a_{i_\star}=0.
    \]
    We may therefore remove $i_\star$ from \eqref{eq:peeling_relation} and apply the private-pole lemma to $J\setminus\{i_\star\}$. Repeating the argument gives
    \begin{equation*}\label{eq:peeling_conclusion}
        \dot q_j^\ell\equiv0 \quad \text{whenever } a_j\neq 0,
        \qquad
        \dot a_j=0,
        \qquad j\in J,
    \end{equation*}
    and finally $c=0$.

    We now apply this argument layer by~layer. Recall \eqref{eq:deep_tanh_layers}. Since
    \[
        D\mathcal R(\theta)[\dot\theta]
        =
        \sum_{i\in[P]} \dot\theta_i 
        \begin{pmatrix}
        \frac{\partial}{\partial\theta_i}u_{\theta,1}\\
        \vdots\\
        \frac{\partial}{\partial\theta_i}u_{\theta,n_{L+1}}
        \end{pmatrix}
        \equiv 0,
    \]
    we have
    \[
        \dot q_i^{L+1}\equiv0,
        \qquad i\in[n_{L+1}].
    \]
    We now descend through the layers, starting with $\ell=L$. Suppose now inductively that
    \[
        \dot q_i^{\ell+1}\equiv0,
        \qquad i\in[n_{\ell+1}].
    \]
    For every such $i$,
    \[
        0
        =
        \dot b_i^{\ell+1}
        +
        \sum_{j\in[n_\ell]}
        \left(
            \dot w_{ij}^{\ell+1}h_j^\ell+
            w_{ij}^{\ell+1}\dot h_j^\ell
        \right).
    \]
    Applying the peeling argument with
    \[
        J=[n_\ell],
        \qquad
        a_j=w_{ij}^{\ell +1}
    \]
    gives
    \[
        \dot w_{ij}^{\ell+1}=0,
        \qquad
        \dot b_i^{\ell+1}=0
    \]
    and
    \[
        \dot q_j^\ell\equiv0
        \qquad
        \text{for every } j \text{ with }  w_{ij}^{\ell+1} \neq 0.
    \]
    By~\eqref{eq:every_neuron_has_outgoing_edge}, for every $j\in[n_\ell]$ there exists $i\in[n_{\ell+1}]$ such that $w_{ij}^{\ell+1}\neq0.$ Hence
    \[
        \dot q_j^\ell\equiv0,
        \qquad
        j\in[n_\ell],
    \]
    so doing this for all $i\in[n_{\ell+1}]$ yields
    \[
        \dot W^{\ell+1}=0,
        \qquad
        \dot b^{\ell+1}=0,
        \qquad
        \dot q_j^\ell\equiv0
        \quad j\in[n_\ell].
    \]

    Descending through the layers gives
    \[
        \dot q_i^1\equiv0,
        \qquad i\in[n_1].
    \]
    Since
    \[
        q_i^1(x)=w_i^1x+b_i^1,
    \]
    we obtain
    \[
        \dot w_i^1=\dot b_i^1=0.
    \]
    Hence every component of~$\dot\theta$ vanishes, i.e.,
    \[
        \dot\theta=0.
    \]
\end{proof}

\begin{proof}[Proof of~Theorem \ref{thm:full_rank_realization_d}]
    We first prove that for $\theta \in \Theta^{\tanh}_{\mathrm{reg}}$, $D\mathcal R(\theta)$ is injective, i.e.,
    \[
        D\mathcal R(\theta)[\dot\theta]\equiv0
        \quad\text{on }\R^d
        \qquad\Longrightarrow\qquad
        \dot\theta=0.
    \]
    Thus, choose $x_0\in\R^d$ such that
    \[
        q_i^1(x_0)\neq\pm q_j^1(x_0)
        \qquad 
        \text{for }
        i\neq j,
    \]
    and choose linearly independent $v_1,\dots,v_d\in\R^d$ such that
    \[
        \langle w_i^1,v_r\rangle\neq0
        \qquad
        \text{for }
        i\in[n_1],\ r\in[d].
    \]
    Such choices exist since regularity implies
    \[ 
    (w_i^1,b_i^1)\neq\pm(w_j^1,b_j^1), \qquad i\neq j, 
    \] 
    hence the excluded sets are finite unions of~proper affine hyperplanes.

    Fix $r\in[d]$ and restrict the network to
    \[
        x=x_0+tv_r,\qquad t\in\R.
    \]
    Its first-layer preactivations are
    \[
        q_i^1(x_0+tv_r)
        =
        \langle w_i^1,v_r\rangle t+q_i^1(x_0).
    \]
    By~the choice of~$x_0$ and $v_r$, the restricted network is a~regular one-input network, so Lemma~\ref{lem:full_rank_realization} applies.

    Under this restriction, the first-layer parameter variations become
    \[
        \langle\dot w_i^1,v_r\rangle,
        \qquad
        \langle\dot w_i^1,x_0\rangle+\dot b_i^1,
    \]
    while all parameter variations in the deeper layers are unchanged. Hence
    \[
        \langle\dot w_i^1,v_r\rangle=0,
        \qquad
        \langle\dot w_i^1,x_0\rangle+\dot b_i^1=0,
    \]
    for every $i\in[n_1]$, and all parameter variations in layers $2,\dots,L+1$ vanish. Applying this for $r=1,\dots,d$ and using that $v_1,\dots,v_d$ form a~basis gives
    \[
        \dot w_i^1=0.
    \]
    Hence also $\dot b_i^1=0$, and therefore
    \[
        \dot\theta=0.
    \]

    For the converse, let $\theta\notin\Theta_{\mathrm{reg}}^{\tanh}$. Recall that $D\mathcal R^{\tanh}(\theta)[\dot\theta]=\dot q^{L+1}$ and $\dot h^0\equiv0$, since $h^0(x)=x$. We construct $\dot\theta\neq0$ for which $\dot q^{L+1}$ vanishes. In each case below, all unspecified components of $\dot\theta$ are set to zero. If $\dot W^m=0$ and $\dot b^m=0$ for $m=1,\dots,\ell$, then \eqref{eq:variation_h_simple} and \eqref{eq:variation_q} give, inductively for $m=1,\dots,\ell$,
    \[
        \dot q_i^m=\sum_{j=1}^{n_{m-1}}w_{ij}^m\dot h_j^{m-1}\equiv0,
        \qquad
        \dot h_i^m=\bigl(1-(h_i^m)^2\bigr)\dot q_i^m\equiv0.
    \]
    Because $\theta$ is not regular, at least one of the following cases occurs for some $\ell\in[L]$:
    
    \begin{enumerate}
        \item If \eqref{eq:every_neuron_has_parent} is violated, there is some $i\in[n_\ell]$ satisfying $W^\ell_{i,:}=0$. Then $h_i^\ell\equiv c:=\tanh(b_i^\ell)$. Choose $v\in\R^{n_{\ell+1}}\setminus\{0\}$ and set
        \[
            \dot W^{\ell+1}_{:,i}=v,
            \qquad
            \dot b^{\ell+1}=-cv.
        \]
        All coefficients in layers $1,\dots,\ell$ vanish, so the preceding induction gives $\dot h^\ell\equiv0$. Substituting the chosen coefficients into \eqref{eq:variation_q}, we obtain
        \[
            \dot q^{\ell+1}=-cv+vh_i^\ell=v(h_i^\ell-c)\equiv0.
        \]
    
        \item If \eqref{eq:every_neuron_has_outgoing_edge} is violated, there is some $i\in[n_\ell]$ satisfying $W^{\ell+1}_{:,i}=0$. Set $\dot b_i^\ell=1$. All coefficients in layers $1,\dots,\ell-1$ vanish, so $\dot h^{\ell-1}\equiv0$. Consequently, \eqref{eq:variation_q} gives
        \[
            \dot q_i^\ell=1,
            \qquad
            \dot q_j^\ell=0 \quad \text{for } j\neq i.
        \]
        By \eqref{eq:variation_h_simple}, it follows that $\dot h_i^\ell=1-(h_i^\ell)^2$ and $\dot h_j^\ell=0$ for $j\neq i$. Since all coefficients in layer $\ell+1$ vanish, we obtain
        \[
            \dot q^{\ell+1}=W^{\ell+1}_{:,i}\bigl(1-(h_i^\ell)^2\bigr)\equiv0.
        \]
    
        \item If \eqref{eq:regularity_of_params_2} is violated, there are distinct $i,j\in[n_\ell]$ and some $\varepsilon\in\{\pm1\}$ satisfying
        \[
            (W^\ell_{i,:},b_i^\ell)=\varepsilon \, (W^\ell_{j,:},b_j^\ell).
        \]
        Then $q_i^\ell=\varepsilon q_j^\ell$ and, by oddness of $\tanh$, $h_i^\ell=\varepsilon h_j^\ell$. Choose $v\in\R^{n_{\ell+1}}\setminus\{0\}$ and set
        \[
            \dot W^{\ell+1}_{:,i}=v,
            \qquad
            \dot W^{\ell+1}_{:,j}=-\varepsilon v.
        \]
        Again, all coefficients in layers $1,\dots,\ell$ vanish, so $\dot h^\ell\equiv0$. Therefore, \eqref{eq:variation_q} yields
        \[
            \dot q^{\ell+1}=vh_i^\ell-\varepsilon vh_j^\ell\equiv0.
        \]
    \end{enumerate}
    
    In every case, $\dot\theta\neq0$ and $\dot q^{\ell+1}\equiv0$. All coefficients in layers $\ell+2,\dots,L+1$ vanish. Thus, for each $m=\ell+2,\dots,L+1$, the identity $\dot q^{m-1}\equiv0$ implies $\dot h^{m-1}\equiv0$ by \eqref{eq:variation_h_simple}, and then $\dot q^m=W^m\dot h^{m-1}\equiv0$ by \eqref{eq:variation_q}. Induction therefore gives
    \[
        D\mathcal R^{\tanh}(\theta)[\dot\theta]=\dot q^{L+1}\equiv0.
    \]
    Thus the functions $\partial_{\theta_j}u_\theta$, $j\in[P]$, are linearly dependent, proving the converse and hence the asserted equality.
\end{proof}
\subsection{Stable sampling}
\begin{lemma}[{\cite[Proposition~2]{vlavcic2021affine}}]\label{lem:global_identifiability_tanh}
    Global identifiability for $(\mathcal{R}^{\tanh},\Theta^{\tanh}_{\mathrm{reg}})$ modulo $\Gamma^{\tanh}$ holds true. 
\end{lemma}
%
%
\begin{proof}[Proof of~Theorem~\ref{thm:stable_sampling_tanh}]
    By~the definition of~regularity, $\Theta_{\mathrm{reg}}^{\tanh}$ is described by~finitely many equalities and inequalities. Hence it is semialgebraic and therefore $\sigma$-subanalytic; see Definition~\ref{def:semi/sub-analyticity}. Moreover, the mappings
    \[
        (\theta,x)\mapsto\mathcal R^{\tanh}(\theta)(x)
        \quad\text{and}\quad
        (\theta,v,x)\mapsto
        D\mathcal R^{\tanh}(\theta)[v](x)
    \]
    are real analytic and hence $\sigma$-subanalytic. In particular, they are real analytic in $x$ on $Y:=\R^d$ and $\nu\mapsto\mathcal R^{\tanh}(\nu)(x)$ is $C^1$ for every $x\in K$. Furthermore, global identifiability of~$(\mathcal R^{\tanh},\Theta_{\mathrm{reg}}^{\tanh})$ modulo $\Gamma^{\tanh}$ follows from Lemma~\ref{lem:global_identifiability_tanh}. Since $\Theta_{\mathrm{reg}}^{\tanh}=\Theta_{\mathrm{reg,fr}}^{\tanh}$ by~Theorem~\ref{thm:full_rank_realization_d}, the assumptions of~Theorem~\ref{thm:generic_global_sampling} are satisifed with parameter set $\Theta_{\mathrm{reg}}^{\tanh}$ and $Y=\R^d$, which yields the claim.
\end{proof}

\section*{Declaration of~AI assistance}
During the preparation of~this manuscript, the authors used large language models (GPT 5.6 Pro and 6.1 Sol) for literature searches, language editing, and assistance with the proof of~Lemma~\ref{lem:private_poles}. The authors supplied a~proof for the shallow case $L=1$ and outlined the main idea for the general case, namely to proceed by~induction on $L$ and use \cite[Proposition~8]{vlavcic2021affine} to construct poles of~nontrivial linear combinations of~tangent vectors. Based on this input, an~LLM generated a~proof sketch for the general case, which the authors developed into the proof presented in the manuscript. The authors reviewed and verified the resulting content and take full responsibility for the final manuscript.

\newpage
\printbibliography

\end{document}